\documentclass[a4paper,12pt]{article}
	\usepackage[english]{babel}
	\usepackage{blindtext}
	\usepackage{latexsym}
	\usepackage{amscd}
	\usepackage{amssymb}
	\usepackage{amstext}
	\usepackage{amsmath}
	\usepackage{amsxtra}
	\usepackage{fancyhdr}
	\usepackage{makeidx}
	\usepackage{multirow}
	\usepackage{multicol}
	\usepackage{graphicx}
	\usepackage{epstopdf}
	\usepackage{ulem}
	\usepackage{hyperref}
	\usepackage{epsfig}
	\usepackage{color}
	\usepackage{enumerate}
	\usepackage[table]{xcolor}
	\usepackage{authblk}
	\newcommand{\etal}{\textit{et al.}}
	\newcommand{\bigO}{\mathcal{O}}
	
	\newcommand{\Var}{\mathbb{V}ar}
	\newcommand{\pr}{\mathbb{P}}
		
	\newcommand\smallO{
		\mathchoice
		{{\scriptstyle\mathcal{O}}}% \displaystyle
		{{\scriptstyle\mathcal{O}}}% \textstyle
		{{\scriptscriptstyle\mathcal{O}}}% \scriptstyle
		{\scalebox{.7}{$\scriptscriptstyle\mathcal{O}$}}%\scriptscriptstyle
	}

	\DeclareMathOperator{\bias}{bias}
	\DeclareMathOperator{\cov}{Cov}
	\newtheorem{pro}{Proposition}[section]
	\newtheorem{theorem}{Theorem}[section]
	\newtheorem{lemma}[pro]{Lemma}
	
	\newtheorem{hq}[pro]{Corollary}
	\newtheorem{ly}[pro]{Remark}
\begin{document}
	\fontsize{12pt}{16pt}\selectfont
	\title{\textbf{Nonparametric density estimation for stationary processes under multiplicative measurement errors}}
	
	\author[a,b]{Dang Duc Trong}
	\author[a,b]{Hoang Van Ha}
    \author[a,b,c]{Thai Phuc Hung}

   \affil[a]{Faculty of Mathematics and Computer Science, University of Science, Ho Chi Minh City, Vietnam.} \affil[b]{Vietnam National University, Ho Chi Minh City, Vietnam.} \affil[c]{Faculty of Basics, Soc Trang Community College, Soc Trang Province, Viet Nam.}
	\maketitle
	
	\begin{abstract}
This paper focuses on estimating the invariant density function $f_X$ of the strongly mixing stationary process $X_t$ in the multiplicative measurement errors model $Y_t = X_t U_t$, where $U_t$ is also a strongly mixing stationary process. We propose a novel approach to handle non-independent data, typical in real-world scenarios. For instance, data collected from various groups may exhibit interdependencies within each group, resembling data generated from $m$-dependent stationary processes, a subset of stationary processes. This study extends the applicability of the model $Y_t = X_t U_t$ to diverse scientific domains dealing with complex dependent data. The paper outlines our estimation techniques, discusses convergence rates, establishes a lower bound on the minimax risk, and demonstrates the asymptotic normality of the estimator for $f_X$ under smooth error distributions. Through examples and simulations, we showcase the efficacy of our estimator. The paper concludes by providing proofs for the presented theoretical results.
	\end{abstract}
	
% new LaTeX commands

%1. INTRODUCTION AAAAAAAAAAAAAAAAAAAAAAAAAAAAAAAAAAAAAAAAAAAAAAAAAA
%1. INTRODUCTION AAAAAAAAAAAAAAAAAAAAAAAAAAAAAAAAAAAAAAAAAAAAAAAAAA

\section{Introduction}

Let $X_{t_1}, \ldots, X_{t_n}$ be positive random variables generated from the strongly mixing stationary process $X_t$ with unknown invariant density function $f_X: \mathbb{R}_+ \to \mathbb{R}_+$. Our objective is to estimate the density function $f_X$ based on observations, while accounting for the presence of multiplicative measurement errors. In particular, let us assume that we have a sample consisting of observations $Y_{t_1}, \ldots, Y_{t_n}$ according to the following model
\begin{equation}
Y_{t_j} = X_{t_j} U_{t_j}, \, j=1, \ldots, n. \nonumber
\end{equation}
Here $U_{t_1}, \ldots, U_{t_n}$ are error random variables generated from the strongly mixing stationary processes $U_t$ with known invariant density function $g: \mathbb{R}_+ \to \mathbb{R}_+$. The processes $X_t$ and $U_t$ are independent. The density function $g$ is called error density. Assume that the strongly mixing stationary processes are observed at discrete time $t_j = j\Delta$ (where $\Delta$ is a positive constant). For brevity, we will denote $X_{t_j}, Y_{t_j}, U_{t_j}$ by $X_j, Y_j, U_j$, respectively. Consequently, our model can be expressed as
\begin{equation}
Y_{j} = X_{j} U_{j}, \, j=1, \ldots, n. \label{eq:11}
\end{equation}
In this setting the invariant density function $f_Y: \mathbb{R}_+ \to \mathbb{R}_+$ of the process $\{Y_j\}$ is given by
$$f_Y(y) = [f_X * g](y) =\int_0^\infty {f_X(x)g(y/x)x^{-1}dx}, \, \forall y \in \mathbb{R}_+ $$
such that $*$ denotes multiplicative convolution. The estimation of $f_X$ using a sample $Y_1, \ldots, Y_n$ generated from the strongly  mixing stationary process $Y_t = X_t U_t$ with density function $f_Y$ is thus an inverse problem called ``multiplicative deconvolution".

Model \eqref{eq:11} has found a lot of  applications in various fields of research and has been extensively studied in the context where $X_1, \ldots, X_n$ are independent and identically distributed (i.i.d.) samples from the distribution of $X$, and $U_1, \ldots, U_n$ are i.i.d. samples from the distribution of $U$. In the field of survival analysis, as elucidated and motivated in Van Es \etal  \cite{Es2000}, when $X$ represents the true positive survival time of a patient who has been sampled, and the random variable $Y$ represents their observed survival time, it is reasonable to assume that the sampling time is uniformly distributed throughout the entire survival period of length $X$. Hence, we have $Y = XU$, where $U$ follows a uniform distribution on the interval $\left[0, 1\right]$. Additionally, in survey sampling, model \eqref{eq:11} is also employed to investigate sensitive information among respondents, where obtaining accurate responses can be challenging, as discussed in Shou and Gupta \cite{Shou2023}.

In the traditional framework where $X_1, \ldots, X_n \mathop \sim \limits^{i.i.d} X$ and $U_1, \ldots, U_n \mathop \sim \limits^{i.i.d} U$, the foundational works of Vardi \cite{Vardi1989} and Vardi \etal \cite{Vardi1992} were pivotal in introducing and investigating deeply into multiplicative censoring, a specific aspect of the multiplicative deconvolution problem where the multiplicative error U is uniformly distributed over the interval $\left[0, 1\right]$. The task of estimating the cumulative distribution function of $X$ was significantly advanced in Vardi \etal \cite{Vardi1992} and Asgharian and Wolfson \cite{Asgharian2005}. Meanwhile, series expansion methods, approaching the model as an inverse problem, were explored in Andersen and Hansen \cite{Andersen2001}. In the realm of density estimation within a multiplicative censoring model, Brunel \etal \cite{Brunel2016}  introduced both kernel estimators and convolution power kernel estimators. The study of Comte and Dion \cite{Comte2016} focuses on a projection density estimator based on the Laguerre basis, assuming a uniform error distribution across the interval $\left[1-\alpha, 1+\alpha\right]$, where $\alpha$ lies within $\left(0, 1\right)$. In the work of Belomestny and Goldenshluger \cite{Belomestny2016}, the authors study the scenario where the error $U$ is beta-distributed. Additionally, Miguel \etal  \cite{Miguel2021} highlighted the transformative role of Mellin transform in a fully data-driven procedure for nonparametric density estimation. Lastly, Miguel \etal  \cite{Miguel2023} focuses on nonparametric estimation using a plug-in approach, combining Mellin transform estimation and regularization for a linear functional evaluated at an unknown density function. 

In the context of dependent data where $X_t$ and $U_t$ are stationary proccesses, the considering model has received limited attention in research studies. However, in the scenario where $X_1, \ldots, X_n$ being a stationary process and $U_1, \ldots, U_n \mathop \sim \limits^{i.i.d} U$, the work by Miguel and Phandoidaen \cite{Miguel2022} stands out, presenting an estimation for the unknown survival function within the framework of model \eqref{eq:11}.

In this study, our primary focus revolves around the estimation of the invariant density function $f_X$ within scenarios where both $X_t$ and $U_t$ are conceptualized as stationary processes. A pivotal aspect of our investigation lies in the adept handling of non-independent data, mirroring real-world complexities. A prevalent scenario arises when data is gathered from diverse groups, and within each group, there exist interdependencies among the data points. These interdependencies align with data patterns generated from $m$-dependent stationary processes, constituting a distinctive subset of stationary processes. The incorporation of this assumption, especially concerning $U_t$, introduces an innovative approach that has not been explored previously, enriching the precision and profundity of our modeling and density estimation. This study not only broadens the applications of model \eqref{eq:11} to various scientific domains but also proves invaluable in instances involving intricate interdependencies among data points. This paper focuses on the smooth error density class, which includes many important distributions such as uniform, beta, etc., under the assumption
\begin{eqnarray}
c{\left( {1 + \left| p \right|} \right)^{ - \kappa }} \le \left| {{g^{mt}}\left( p\right)} \right| \le C{\left( {1 + \left| p \right|} \right)^{ - \kappa }}, \label{eq:12}
\end{eqnarray}
where $C \ge c > 0, \kappa > 0$ and ${g^{mt}}\left( p \right)$ is the Mellin transform of the invariant error density function $g$, see for instance Miguel \etal \cite{Miguel2021}.

The structure of this paper is organized as follows. Section 2 outlines our proposed estimator by using the Mellin transform. In Section 3, we discuss convergence rates, provide a lower bound on the minimax risk, and explore the asymptotic normality of the estimator for $f_X$, all without including proofs. Applications and simulations are thoroughly examined in Section 4. Concluding the paper, Section 5 offers detailed proofs for the results introduced in Section 3.

% 2. ESTIMATOR BBBBBBBBBBBBBBBBBBB
% 2. ESTIMATOR BBBBBBBBBBBBBBBBBBB

\section{Construction of the estimators}

\subsection{Mellin transform}
In the multiplicative censoring model \eqref{eq:11}, we consider a density estimator using the Mellin transform. The key to the analysis of the multiplicative deconvolution problem is the convolution theorem of the Mellin transform $\mathcal{M}$, which roughly states $\mathcal{M}[f_Y] = \mathcal{M}[f_X] \mathcal{M}[g]$ for a density $f_Y = f_X *g$. For $\varrho \ge 1$, we denote by $L_\varrho^+ = L_\varrho(\mathbb{R}_+)$ the set of Lebesgue measurable functions $f$ satisfying $||f||_\varrho=\left(\int_0^\infty |f\left(x\right)|^\varrho dx \right)^{1/\varrho} < \infty$. For $f \in L_1^+$, the Mellin transform of $f$ in the point $c + ip \in \mathbb{C}$ is defined by
\begin{eqnarray}
\mathcal{M}_c[f] \left( p \right):= \mathcal{M}[f] \left( c+ ip \right) = \int_0^{\infty} x^{c-1+ip}f \left( x \right)dx \label{eq:21}
\end{eqnarray}
provided that the integral is absolutely convergent. If there exists a $c \in \mathbb{R}$ such that the mapping $x \mapsto x^{c-1} f(x)$ is integrable over $\mathbb{R}_+$ then the region $\Xi_f \subseteq \mathbb{C}$ of absolute convergence of the integral in \eqref{eq:21} is either a vertical
strip $\{s + ip \in \mathbb{C} : s \in (a, b), p \in \mathbb{R}\}$ for $a < b$ with $c \in (a, b)$ or a vertical line $\{c + ip \in \mathbb{C} : p \in \mathbb{R}\}$. We can see some techniques to determine $\Xi_f$ in (\cite{Miguel2021}). Note that for any density $f \in L_1^+$ the vertical line $\{1 + ip \in \mathbb{C} : p \in \mathbb{R}\}$ belongs to $\Xi_f$, and hence the Mellin transform $\mathcal{M}_1[f]$ is well-defined. In this article, we will restrict to the case of model parameter $c = 1$ and denote $f^{mt} = \mathcal{M}_1[f]$. Nevertheless, we can extend for arbitrary model parameters $c \in \mathbb{R}$. For $c = 1$, the inversion formula of the Mellin transform is given by
\begin{equation}
 f\left( x \right) = \frac{1}{{2\pi }}\int_{ - \infty }^{ + \infty } {{x^{- 1 - ip}} f^{mt} \left( p \right) dp}, \quad x \in \mathbb{R}_+. 
\label{eq:22}
\end{equation} 

\subsection{The proposed estimators}

In the present paper, we consider the problem on a dependent process. As known, there are various types of dependence conditions which is deeply studied (see Bradley \cite{Bradley2005}).  In this paper, we limit ourselves to the $\alpha-$mixing dependence which was introduced by Rosenblatt 
\cite{Rosenblatt1956}. Let $\mathcal{F}_i^k$ be the $\sigma$-algebra of events generated by the random variables $\left\{ X_j, 1\leq i \le j \le k \right\}$. The process $\left\{ X_j \right\}$ is called strongly mixing or $\alpha$-mixing (see \cite{Bradley2005, Rosenblatt1956}) if
$$\sup_{j\in\mathbb{N}}\mathop {\sup }\limits_{\scriptstyle A \in \mathcal{F}_{1 }^j\hfill\atop \scriptstyle B \in \mathcal{F}_{j+k}^{ + \infty }\hfill} \left| {\mathbb{P}\left[ {AB} \right] - \mathbb{P}\left[ A \right].\mathbb{P}\left[ B \right]} \right| = \alpha_X \left( k \right) \to 0 \text{ as } k \to  + \infty, $$
where $\alpha_X \left( k \right)$ is the strong mixing coefficient of the process $\left\{ X_j \right\}$.
The quantity is a kind of measures of
dependence.  When $(X_j)$ is mutually independent then $\sup_k\alpha_X \left( k \right)=0$. The condition 
$\alpha_X(k)\to 0$ means that  $X_j$ and $X_{j+k}$ are ``asymptotically independent'' as 
$k\to\infty$. With the dependence properties, we can establish an estimator of deconvolution which is consistent. 

From \eqref{eq:22}, to obtain an estimator of $f_X(x)$, we first try to construct a suitable estimator of $f_X^{mt}(p) = f_Y^{mt}(p)/g^{mt}(p)$. By replacing $f_Y^{mt}(p)$ with its empirical counterpart $\widehat{f}_Y^{mt}(p) = \frac{1}{n} \sum_{j=1}^n Y_j^{ip}$ and using the kernel function construction as in Belomestny and Goldenshluger \cite{Belomestny2020}, we obtain an estimator of $f_X(x)$ in the final form 
\begin{equation}
 {\widehat f}_{n} \left(x \right) = \frac{1}{{2\pi }} \int_{ - \infty }^{ + \infty } {{x^{ - 1 - ip}}\frac{K^{ft} \left(pb_n \right)}{g^{mt}(p)}\frac{1}{n}\sum\limits_{j = 1}^n {{Y_j^{ip}}}} dp, \quad x \in \mathbb{R}_+,
\label{eq:23}
\end{equation}
where $b_n$ is the positive parameter which tends to $0$ as $n$ tends to infinity, $K$ is a known kernel function and $K^{ft}\left(p\right) := \int_{-\infty}^{+\infty}e^{ipx}K\left(x\right)dx$ is the Fourier transform of $K$. 

\subsection{Assumptions}

\noindent \textbf{Assumption A}
\begin{enumerate}
\item[$\left( A_1\right)$] The kernel function $K$ is supported on $\left[-1, 1\right]$, bounded and for a positive integer number $m$,
\begin{eqnarray}
\int_{-1}^1 K\left(x\right)dx = 1, \quad \int_{-1}^1x^kK\left(x\right)dx = 0, \quad k=1,...,m. \nonumber
\end{eqnarray}
\item[$\left( A_2\right)$] The Fourier transform $K^{ft}$ of the kernel function $K$ satisfies $\int{|p|^{k \kappa}\left|K^{ft}\left(p\right)\right|^k}dp < \infty $ for $k=1,2$ and $\kappa$ is given in \eqref{eq:12}.
\end{enumerate}
 
\noindent \textbf{Assumption B}
\begin{enumerate}
\item[$\left( B_1\right)$] For some $\kappa > 0, 0 < c \le C, D > 0$, the error density function $g$ of the process $U_t$ has the Mellin transform satisfying $$c{\left( {1 + \left| p \right|} \right)^{ - \kappa }} \le \left| {{g^{mt}}\left( p\right)} \right| \le C{\left( {1 + \left| p \right|} \right)^{ - \kappa }}, \left| {{{\left( {{g^{mt}}\left( p \right)} \right)}^\prime }} \right| \le {D}{\left( {1 + \left| p \right|} \right)^{ - \kappa }} \text{ for all } p \in \mathbb{R}.$$
\item[$\left( B_2\right)$] The 2-dimensional probability density function ${f_{{Y_{j}},{Y_{k}}}}\left( {u,v} \right)$ exists and is bounded for all $1 \le j,k \le n.$
\end{enumerate}

\begin{ly}
	The kernel function $K$ can be constructed satisfying the Assumption A in different ways. For example, let $w\left(x\right) = e^{-x^2/2}/\sqrt{2\pi}$ and for a fixed natural number $m$ let
	\begin{eqnarray}
		K\left(x\right) = \sum\limits_{j = 1}^{m + 1} {\left( {\begin{array}{*{20}{c}}
					{m + 1}\\
					j
			\end{array}} \right){{\left( { - 1} \right)}^{j + 1}}\frac{1}{j}w\left( {\frac{x}{j}} \right)}. \label{eq:24}
	\end{eqnarray}
	The Fourier transform of $K$ is given by
	\begin{eqnarray}
		K^{ft}\left(p\right) = \sum\limits_{j = 1}^{m + 1} {\left( {\begin{array}{*{20}{c}}
					{m + 1}\\
					j
			\end{array}} \right){{\left( { - 1} \right)}^{j + 1}} e^{-j^2p^2/2} }. \label{eq:25}
	\end{eqnarray}
	It is well-known that the kernel $K$ defines in \eqref{eq:24} and its Fourier transform in \eqref{eq:25} satisfying the Assumption A (see Kerkyacharian \etal \cite{Kerkyacharian2001}).
\end{ly}

% 3. PRELIMINARY RESULTS CCCCCCCCCCCCCCCCCCCCCCCCCCCCCCCCCCCCCC
% 3. PRELIMINARY RESULTS CCCCCCCCCCCCCCCCCCCCCCCCCCCCCCCCCCCCCC

\section{Preliminary results}

%\subsection{Some results on the bias and variance of the estimators}
In this Section, we introduce some results on the bias and variance, then the almost surely  converrgence of our proposed estimator. 

The bias of the estimators will be analyzed under a local smoothness assumption on $f_X$. Let $s>0, A>0, x>0$ and $r>1$. We say that $f \in \mathcal{F}_{x,r}\left(A,s\right)$ if $f$ is a probability density, that is, $l=\left\lfloor s \right\rfloor := \max\left\{k \in \mathbb{N}_0:k<s\right\}$ times continuously differentiable, and $\max_{k=1,...,l}\left|f^{(k)}\left(x_1\right)\right| \le A$,
\begin{eqnarray}
\left|f^{(l)}\left(x_1\right) - f^{(l)}\left(x_2\right)\right| \le A\left|x_1 - x_2\right|^{s-l}, \quad \forall x_1,x_2 \in \left[r^{-1}x, rx\right]. \nonumber
\end{eqnarray}

The following Theorem gives an upper bound on the bias of $\widehat{f}_n\left(x\right)$ when $f_X\left(x\right) \in \mathcal{F}_{x,r}\left(A,s\right)$. We refer Lemma 2 and Lemma 3 of \cite{{Belomestny2020}} for the proof of this theorem. 

%The proof of this Theorem can be seen in (\cite{{Belomestny2020}}, Lemma 2 and Lemma 3).

%Theorem 3.1 1111111111111111111111111111111111111111111111111111111111
%Theorem 3.1 1111111111111111111111111111111111111111111111111111111111

\begin{theorem}\label{th:bias}
Let $f_X\left(x\right) \in \mathcal{F}_{x,r}\left(A,s\right)$ and $\widehat{f}_{n}(x)$ be given in \eqref{eq:23}. Let $K$ be the kernel function $K$ satisfying Assumption A with $m \ge \left\lfloor s\right\rfloor + 1$, then for any $x > 0$ and $n$ large enough
\[ \sup\limits_{f_X \in \mathcal{F}_{x,r}\left(A,s\right)} \bias\left[ {{\widehat f }_{n}(x)} \right] \le \bigO\left(b_n^s\right), \]
where $\bias \left[ {{\widehat f }_{n}(x)} \right] = \left| \mathbb{E}{{\widehat f }_{n}(x)} - f_X (x) \right| $.
\end{theorem}

The general bounds for the variance of $\widehat{f}_n$ are given in the Theorem \ref{th:variance}.

% THEOREM 3.2 2222222222222222222222222222222222222222222222222222
% THEOREM 3.2 2222222222222222222222222222222222222222222222222222

\begin{theorem} \label{th:variance}
	Let Assumption A and Assumption B hold. For $\delta > 0$, assume that the processes $X_t, U_t$ have strong mixing coefficients ${\alpha_X }\left( k \right),$ ${\alpha_U }\left( k \right) \le \bigO\left( {{1/{{k^{2 + \delta }}}}} \right)$. If $f_X\left(x\right) \in \mathcal{F}_{x,r}\left(A,s\right)$, then for any positive sequence $\left( {{b_n}} \right)$ such that ${b_n} \to 0$ and $nb_n^{1 + 2\kappa} \to \infty $ as $n \to \infty $,
\begin{eqnarray}
\sigma_1^2 \left(x\right) = \mathop {\lim }\limits_{m \to \infty } \mathop {\inf }\limits_{n \ge m} \left[ {nb_n^{1 + 2\kappa} \mathop {\rm var } \left( {{\widehat f}_n}\left(x\right) \right)} \right] \le \mathop {\lim }\limits_{m \to \infty } \mathop {\sup }\limits_{n \ge m} \left[ {nb_n^{1 + 2\kappa} \mathop {\rm var } \left( {{\widehat f}_n}\left(x\right) \right)} \right] = \sigma_2^2 \left(x\right) \nonumber
\end{eqnarray}
at continuity points $x>0$ of $f_Y$. Where 
\begin{eqnarray}
\sigma_1^2 \left(x\right) & = & \frac{f_Y\left( x \right)}{2\pi x} \mathop {\lim }\limits_{m \to \infty } \mathop {\inf }\limits_{n \ge m} \int\limits {{{\left| {\frac{{b_n^\kappa K^{ft}\left( p \right)}}{{\overline {{g^{mt}}\left( {{p/{{b_n}}}} \right)} }}} \right|}^2}dp} \ge \frac{{f_Y\left( x \right)}}{2\pi x{C_2^2}} \int{|p|^{2\kappa}\left|K^{ft}\left(p\right)\right|^2}dp , \label{eq:31} \\
\sigma_2^2 \left(x\right) & = & \frac{f_Y\left( x \right)}{2\pi x} \mathop {\lim }\limits_{m \to \infty } \mathop {\sup }\limits_{n \ge m} \int\limits {{{\left| {\frac{{b_n^\kappa K^{ft}\left( p \right)}}{{\overline {{g^{mt}}\left( {{p/{{b_n}}}} \right)} }}} \right|}^2}dp} \le \frac{{f_Y\left( x \right)}}{2\pi x{C_1^2}}\int{|p|^{2\kappa}\left|K^{ft}\left(p\right)\right|^2}dp \label{eq:32}
\end{eqnarray}
and $$0 < {C_1} = \mathop {\lim }\limits_{u \to \infty } \mathop {\inf }\limits_{\left| p \right| \ge u} \left| {{g^{mt}}\left( p \right)} \right|{\left| p \right|^\kappa } \le \mathop {\lim }\limits_{u \to \infty } \mathop {\sup }\limits_{\left| p \right| \ge u} \left| {{g^{mt}}\left( p \right)} \right|{\left| p \right|^\kappa }  = {C_2}.$$
\end{theorem}

%\subsection{Almost surely convergence}

In Theorem \ref{th: converge-a.s-estimators}, we show that the estimator $\widehat{f}_n$ converge almost surely to $f_X$ when $n \to \infty$ under some assumptions.

% THEOREM 3.3 3333333333333333333333333333333333333333333333333333
% THEOREM 3.3 3333333333333333333333333333333333333333333333333333

\begin{theorem}\label{th: converge-a.s-estimators}
Let the Assumption A and assumption $B_1$ hold. For $a > 0, 0 < \rho  < 1 $, assume that the stationary processes $X_t, U_t$ have strong mixing coefficients ${\alpha_X }\left( k \right),$ ${\alpha_U }\left( k \right) \le a{\rho ^k}$. If $f_X\left(x\right) \in \mathcal{F}_{x,r}\left(A,s\right)$, then for any positive sequence $\left( {{b_n}} \right)$ such that ${b_n} \to 0$ and $nb_n^{1 + 2\kappa } \to  + \infty $ as $n \to  + \infty $, ${\widehat f_n} \left(x \right) \to f_X \left(x \right)$ almost surely.
\end{theorem}

% 4. MAIN RESULTS DDDDDDDDDDDDDDDDDDDDDDDDDDDDDDDDDDDDDD
% 4. MAIN RESULTS DDDDDDDDDDDDDDDDDDDDDDDDDDDDDDDDDDDDDD

\section{Main results}

\subsection{Asymptotic normality}\label{sec:asymptotic-normality}

The result of asymptotic normality for ${\widehat{f}_n}$ is as follows.
% THEOREM 4.1 4444444444444444444444444444444444444444444444444444
% THEOREM 4.1 4444444444444444444444444444444444444444444444444444
\begin{theorem}\label{th:asymptotic-normality-estimators}
Let Assumption A and Assumption B hold. For $\delta > 0$, assume that the processes $X_t, U_t$ have strong mixing coefficients ${\alpha_X }\left( k \right),$ ${\alpha_U }\left( k \right) \le \bigO\left( {{1/{{k^{2 + \delta }}}}} \right)$. Choose the bandwidth ${b_n}$ so that $nb_n^{1 + 2\kappa } \to \infty $ and a sequence of positive integers $\left\{ {{s_n}} \right\}$ such that ${s_n} \to \infty $, ${s_n} = \smallO \left( {{{\left( {n{b_n}} \right)}^{{1/2}}}} \right)$ 
 and ${\left( {n/{{b_n}}} \right)^{1/2}}\alpha_U \left( {{s_n}} \right) \to 0\,\,{\rm{as}}\,\,n \to \infty$. Then 
\begin{eqnarray}
\mathop {\lim }\limits_{m \to \infty } \mathop {\inf }\limits_{n \ge m} \left\{ {{n^{1 \mathord{\left/
 {\vphantom {1 2}} \right.
 \kern-\nulldelimiterspace} 2}}b_n^{{1 \mathord{\left/
 {\vphantom {1 2}} \right.
 \kern-\nulldelimiterspace} 2} + \kappa }\left( {{{\widehat f}_n}\left( x \right) - \mathbb{E}{{\widehat f}_n}\left( x \right)} \right)} \right\}\mathop  = \limits^D \mathcal{N}\left( {0;\sigma _1^2\left( x \right)} \right), \label{eq:33} \\
\mathop {\lim }\limits_{m \to \infty } \mathop {\sup }\limits_{n \ge m} \left\{ {{n^{1 \mathord{\left/
 {\vphantom {1 2}} \right.
 \kern-\nulldelimiterspace} 2}}b_n^{{1 \mathord{\left/
 {\vphantom {1 2}} \right.
 \kern-\nulldelimiterspace} 2} + \kappa }\left( {{{\widehat f}_n}\left( x \right) - \mathbb{E}{{\widehat f}_n}\left( x \right)} \right)} \right\}\mathop  = \limits^D \mathcal{N}\left( {0;\sigma _2^2\left( x \right)} \right), \,\, \label{eq:34}
\end{eqnarray}
at continuity points $x > 0$ of $f_Y$, where $\sigma _1^2\left( x \right)$ and $\sigma _2^2\left( x \right)$ given in \eqref{eq:31} and \eqref{eq:32}.
\end{theorem}

\begin{ly}
	The selection of the bandwidth $b_n$ in Theorem \ref{th:asymptotic-normality-estimators} determines the existence of the sequence of positive integers $\{ s_n \}$. If the stationary processes $X_t, U_t$ have strongly mixing coefficients ${\alpha_X }\left( k \right) = \bigO\left( {{1/{{k^{2 + \delta }}}}} \right)$ or ${\alpha_U }\left( k \right) = \bigO\left( {{1/{{k^{2 + \delta }}}}} \right)$, we can choose the bandwidth $b_n = C n^{-\omega_1} \left( \log n \right)^{\omega_2}$ (where $C>0, \omega_1 < \min \{1/\left( 1+2\kappa \right), \left( 1+\delta \right)/\left( 3+\delta \right) \}, \omega_2 \in \mathbb{R}$ or $C>0, \omega_1 = \min \{1/\left( 1+2\kappa \right), \left( 1+\delta \right)/\left( 3+\delta \right) \},$ $\omega_2 > 0$). In the case ${\alpha_X }\left( k \right),$ ${\alpha_U }\left( k \right) \le \bigO \left(\rho^k\right) $ ($0 < \rho < 1$), bandwidth is chosen satisfying $b_n = C n^{-\omega_1} \left( \log n \right)^{\omega_2}$ ($C>0, \omega_1 < 1/\left( 1+2\kappa \right), \omega_2 \in \mathbb{R}$ or $C>0, \omega_1 = 1/\left( 1+2\kappa \right),\omega_2 > 0$).
\end{ly}

\begin{hq}\label{coro:asymptotic-normality}
Let Assumptions of Theorem \ref{th:asymptotic-normality-estimators} hold. If $f_X\left(x\right) \in \mathcal{F}_{x,r}\left(A,s\right)$, choose the bandwidth ${b_n}$ so that $nb_n^{1 + 2\kappa } \to \infty $ and $nb_n^{1 + 2\kappa + 2s } \to 0$ as $n \to \infty $. Then
\begin{eqnarray}
\mathop {\lim }\limits_{m \to \infty } \mathop {\inf }\limits_{n \geqslant m} \left\{ {{n^{{1/2}}}b_n^{{1/2} + \kappa }\left( {{{\widehat f}_n}\left( x \right) - {f_X }\left( x \right)} \right)} \right\}\mathop  = \limits^D \mathcal{N}\left( {0;\sigma _1^2\left( x \right)} \right), \nonumber \\
\mathop {\lim }\limits_{m \to \infty } \mathop {\sup }\limits_{n \geqslant m} \left\{ {{n^{{1/2}}}b_n^{{1/2} + \kappa }\left( {{{\widehat f}_n}\left( x \right) - {f_X }\left( x \right)} \right)} \right\}\mathop  = \limits^D \mathcal{N}\left( {0;\sigma _2^2\left( x \right)} \right) \,\, \nonumber
\end{eqnarray}
at continuity points $x > 0$ of $f_Y$, where $\sigma _1^2\left( x \right)$ and $\sigma _2^2\left( x \right)$ given in \eqref{eq:31} and \eqref{eq:32}.
\end{hq}

\subsection{Bounds of the estimators}\label{sec:upper-lower-bounds}

The following theorem gives us upper bounds of mean square error of $\widehat{f}_n \left(x\right)$.

% THEOREM 4.2 5555555555555555555555555555555555555555555555555555
% THEOREM 4.2 5555555555555555555555555555555555555555555555555555

\begin{theorem}\label{th:upper-bounds-estimators}
Let Assumption A and Assumption B hold, and assume that $f_X\left(x\right) \in \mathcal{F}_{x,r}\left(A,s\right)$. By choosing the bandwidth ${b_n} = \bigO \left({n^{-1/\left(1 + 2\kappa +2s \right)}}\right)$, we have
\begin{equation}
\sup\limits_{f_X \in \mathcal{F}_{x,r}\left(A,s\right)} \mathbb{E}\left|\widehat {f}_{n} \left(x\right) - f_X \left(x\right) \right|^2 \le \bigO \left( n^{-2s/\left(1 + 2\kappa + 2s \right)} \right), \nonumber
\end{equation}
for $n$ sufficiently large.
\end{theorem}

Let $\mathcal{T}_\alpha$ be the set of all arbitrary estimator $\widehat{f}$ of $f_X$ based on data $Y_1,...,Y_n$ which is generated from the $\alpha-mixing$ stationary process $Y_t$ in model \eqref{eq:11}. The next Theorem gives lower bounds for estimating densities of the set $\mathcal{T}_\alpha$. 

% THEOREM 4.3 6666666666666666666666666666666666666666666666666666
% THEOREM 4.3 6666666666666666666666666666666666666666666666666666

\begin{theorem}\label{th:lower-bounds-estimators}
Let $x \ge B > 0$ for some constant $B$ and suppose that assumption $B_1$ hold and $\kappa > 1/2$. Then
\begin{equation}
\inf\limits_{\widehat{f} \in \mathcal{T}_\alpha} \sup\limits_{f_X \in \mathcal{F}_{x,r}\left(A,s\right)} \mathbb{E}\left|\widehat {f} \left(x\right) - f_X \left(x\right) \right|^2 \ge \bigO \left( n^{-2s/\left(1 + 2\kappa + 2s \right)} \right), \nonumber
\end{equation}
for $n$ sufficiently large.
\end{theorem}
Theorem \ref{th:upper-bounds-estimators} and Theorem \ref{th:lower-bounds-estimators} show that the estimator $\widehat{f}_n$ attains the minimax rate of convergence $n^{-2s/\left(1 + 2\kappa + 2s \right)}$ when the error density $g$ is ordinary smooth with $\kappa > 1/2$.

% 4. EXAMPLES AND SIMULATIONS   DDDDDDDDDDDDDDDDDDDDDDDDDDDDDDDDDDDDDD
% 4. EXAMPLES AND SIMULATIONS   DDDDDDDDDDDDDDDDDDDDDDDDDDDDDDDDDDDDDD

\section{Examples and Simulations}

\subsection{Cox-Ingersoll-Ross process}

{\raggedright The Cox-Ingersoll-Ross process is the solution to the stochastic differential equation}
$$d{X_t} = \left( {\theta_1  - \theta_2 {X_t}} \right)dt + \theta_3 \sqrt {{X_t}} d{W_t}.$$
Under the hypothesis $2\theta_1 > \theta_3^2$, there exists a positive $\tau$ such that $\alpha_X \left(k\right) \le e^{-\tau k}/4$ (see Corollary 2.1 of Genon-Catalot \etal \cite{Genon2000}). The Cox-Ingersoll-Ross process is then an exponentially strong mixing stationary process and its invariant distribution follows a gamma density function 
$${f_X }\left( x \right) = \frac{{{b^{a }}}}{{\Gamma \left( {a } \right)}}{x^{a - 1}}\exp \left( { - b x} \right),$$
where the shape parameter $a = 2\theta_1/\theta_3^2$ and the rate parameter $b = 2\theta_2/\theta_3^2$.

Assume that $X_t$ is exponentially strongly mixing stationary Cox-Ingersoll-Ross processes with $\theta_1 = 1, \theta_2 = 0.5, \theta_3 = 1 \left(a = 2, b = 1\right)$ and noise random variables $U_j$ follow either an uniform distribution on $\left[a_1, b_1\right]$ (denote as $U_{(a_1, b_1)}$) or Beta distribution (denote as $B_{(a_2, b_2)}$).

We now present simulations with different parameter values using the language R. First, we generate sample paths using the "sde" package (see Iacus \cite{Iacus2008}), and then we compute the estimator's values. We choose the kernel function $K$ as defined in (\ref{eq:24}) with its Fourier transform given in (\ref{eq:25}) for the case of $m = 2$. We also select the bandwidth $b_n = n^{-1/\left(1+2\kappa+2s\right)}$ according to Theorem \ref{th:asymptotic-normality-estimators}. The simulations are based on 2000 noisy observations $Y_j = X_jU_j$ of the Cox-Ingersoll-Ross process $X_t$, with 50 replications. The results of the mean square error (MSE) at some points $x > 0$ for the case where $X_t$ is a Cox-Ingersoll-Ross process and $U_t$ is selected from $\left\{ U_{(0,1)}, U_{(0.5, 1.5)}, B_{(1, 2)}, B_{(2, 1)}, B_{(2, 2)} \right\}$ are presented in Table \ref{table1}.

\begin{table}[!htbp]
\caption{The empirical MSE of estimator at some points $x > 0$ in the case $X_t$ is Cox-Ingersoll-Ross process and $U_t \in \left\{ U_{(0,1)}, U_{(0.5, 1.5)}, B_{(1, 2)}, B_{(2, 1)}, B_{(2, 2)} \right\}$.}
\begin{center}
{\renewcommand\arraystretch{1.2}
{\tabcolsep = 0.3cm
\begin{tabular}{ c  c  c  c  c  c }
\noalign{\hrule height 1.5pt}
%\rowcolor[rgb]{0.8,0.8,0.8}
$x$ & $U_{(0,1)}$ & $U_{(0.5, 1.5)}$ & $B_{(1, 2)}$ & $B_{(2, 1)}$ & $B_{(2, 2)}$ \\

%\rowcolor{gray!20}
\hline
 0.5 & 0.02000825   & 0.01151807   & 0.01535441    & 0.01260686   & 0.01174962 \\
 1.0 & 0.01147437   & 0.009493284  & 0.01464701    & 0.009180121  & 0.01356236 \\
 1.5 & 0.001054993  & 0.0007115975 & 0.001157075   & 0.0008065801 & 0.0009686502 \\
 2.0 & 0.001534092  & 0.001253764  & 0.0006557094  & 0.001036627  & 0.0006276501 \\
 2.5 & 0.001789301  & 0.001902348  & 0.001516154   & 0.001864769  & 0.001463214 \\
 3.0 & 0.001974052  & 0.001711986  & 0.001340604   & 0.001830764  & 0.001237362 \\
 3.5 & 0.001127057  & 0.001042727  & 0.0008609882  & 0.001050847  & 0.0008042085 \\
 4.0 & 0.0005810139 & 0.0004889677 & 0.0005396105  & 0.0005298024 & 0.0004642698 \\
 4.5 & 0.0002820801 & 0.0001894068 & 0.0002175196 & 0.0002062239 & 0.0001943284 \\
\noalign{\hrule height 1.5pt}
\end{tabular}}}
\end{center}
\label{table1}
\end{table}

\subsection{The $m$-Dependent stationary process}

A very important kind of dependence considering distance as a measure of dependence, is the $m$-dependence case. A (discrete time) stationary process $\left\{X_j\right\}_{j \in \mathbb{Z}}$ is call $m$-dependent stationary process if two sets of random variables $\left\{...,X_{k-1},X_k\right\}$ and $\left\{X_{h},X_{h+1},...\right\}$ are independent whenever the time gap between them, represented as $h-k$, exceeds the value of $m$. Clearly, $\alpha_X \left(l\right) = 0$ when $l>m$ and the $m$-dependent stationary process is a special case of strongly mixing stationary process.

In this specific example, we assume that $\left\{X_j\right\}_{j \in \mathbb{Z}}$ is a $30$-dependent stationary process and its invariant distribution is Weibull distribution with the shape parameter $a = 2$ and the scale parameter $b = 5$. We will use the language R to generate data for the 30-dependent process $\left\{ X_j\right\}_{j \in \mathbb{Z}}$ as follows. First, we generate 30 independent data points, all following a Weibull distribution with the shape parameter $a= 2$ and the scale parameter $b = 5$. The $31$st data point will either be one of the first 30 data points with a probability of 0.5, or it will be an independent Weibull data point with the same parameters with a probability of 0.5. Starting from the $(30+l)$th data point (where $l \ge 2$), it will select one of the data points from $l$th to $(30+l-1)$th with a probability of 0.5 (if the data point $(30+l-1)$th is the same as data point $(l-1)$th, it will be removed from the list), and it will also have a 0.5 probability of being an independent Weibull data point with the same parameters. The noise random variables $U_j$ are chosen in $\left\{ U_{(0,1)}, U_{(0.5, 1.5)}, B_{(1, 2)}, B_{(2, 1)}, B_{(2, 2)} \right\}$.

With consistent bandwidth and kernel function, as outlined in Section \ref{sec:asymptotic-normality}, the simulations are based on 2000 noisy observations, denoted as $Y_j = X_jU_j$, originating from the $\left\{X_j\right\}_{j \in \mathbb{Z}}$ processes, and this is repeated 50 times. The results of empirical MSE at various points $x > 0$ are given in Table \ref{table2}.

\begin{table}[!htbp]
\caption{The empirical MSE of estimator at various points $x > 0$ in the case $X_t$ is 30-dependent stationary process and $U_t \in \left\{ U_{(0,1)}, U_{(0.5, 1.5)}, B_{(1, 2)}, B_{(2, 1)}, B_{(2, 2)} \right\}$.}
\begin{center}
{\renewcommand\arraystretch{1.2}
{\tabcolsep = 0.3cm
\begin{tabular}{ c  c  c  c  c  c }
\noalign{\hrule height 1.5pt}
%\rowcolor[rgb]{0.8,0.8,0.8}
$x$ & $U_{(0,1)}$ & $U_{(0.5, 1.5)}$ & $B_{(1, 2)}$ & $B_{(2, 1)}$ & $B_{(2, 2)}$ \\

%\rowcolor{gray!20}
\hline
 0.5 & 0.006237416   & 0.002502395   & 0.003739801   & 0.002115358   & 0.001116817 \\
 1.0 & 0.002888191   & 0.001106438   & 0.009513839   & 0.001234231   & 0.009458161 \\
 2.0 & 0.001482346   & 0.001076423   & 0.001115319   & 0.0009224827  & 0.001071403 \\
 3.0 & 0.0006854486  & 0.0003894643  & 0.001294102   & 0.0002997344  & 0.001251179 \\
 4.0 & 0.0005290895  & 0.0004286804  & 0.0001148088  & 0.0003080888  & 0.0001467156 \\
 5.0 & 0.0001766462  & 0.0001377783  & 0.00008996057 & 0.0001675634  & 0.00008717821 \\
 6.0 & 0.00009856583 & 0.00006615849 & 0.0001096502  & 0.0000772704  & 0.00008234754 \\
 7.0 & 0.00004252327 & 0.00003739727 & 0.00002887898 & 0.00002122861 & 0.00002844099 \\
 8.0 & 0.00002642288 & 0.00002186801 & 0.0000406609  & 0.00002247162 & 0.00004048739 \\
\noalign{\hrule height 1.5pt}
\end{tabular}}}
\end{center}
\label{table2}
\end{table}

%%%%%%%%%%%%%%%%%%%%%%%%%%%%%%%%%%%%%%%%%%%%%%%%%%%%%%%%%%%%%%%%%%%%%%%%%%%%%%%%%%%%%%%%%%%%%%%%%%%%%%%%%%%%%%%%%%%%%%%%%%%%%%%%%%%%%%%%%%%%%%%%%%%%%%%%%%%%%%%%%%%%%%%%%%%%%%%%%

% 5. PROOFS   EEEEEEEEEEEEEEEEEEEEEEEEEEEEEEEEEEEEEEEEEEEEEEEEEEEEEEEEEE
% 5. PROOFS   EEEEEEEEEEEEEEEEEEEEEEEEEEEEEEEEEEEEEEEEEEEEEEEEEEEEEEEEEE

\section{Proofs}

\subsection{Proof of Theorem \ref{th:variance}}
Before proving Theorem Theorem \ref{th:variance}, we prove the following technical lemma.

% LEMMA 1
% LEMMA 1

\begin{lemma}\label{lem1}
Asumme that Assumption A and Assumption $B_1$ hold. We have
\begin{eqnarray}
\left( a \right)\,\, \mathop {\lim }\limits_{m \to \infty } \mathop {\inf }\limits_{n \ge m} \frac{1}{{{b_n}}}\int\limits {{\phi _n}\left( {\frac{{\ln x - \ln y}}{{{b_n}}}} \right) f_Y \left( y \right)dy} 
& = & 2\pi x f_Y\left( x \right)\mathop {\lim }\limits_{m \to \infty } \mathop {\inf }\limits_{n \ge m} \int\limits {{{\left| {\frac{{b_n^\kappa K^{ft}\left( p \right)}}{{\overline {{g^{mt}}\left( {{p \mathord{\left/
 {\vphantom {t {{b_n}}}} \right.
 \kern-\nulldelimiterspace} {{b_n}}}} \right)} }}} \right|}^2}dp} \nonumber \\
& \ge & \frac{{2\pi }}{C_2^2} x f_Y\left( x \right)\int{|p|^{2\kappa}\left| K^{ft}(p) \right|^2}dp, \nonumber \\
\left( b \right) \,\, \displaystyle \mathop {\lim }\limits_{m \to \infty } \mathop {\sup }\limits_{n \ge m} \frac{1}{{{b_n}}}\int\limits {{\phi _n}\left( {\frac{{\ln x - \ln y}}{{{b_n}}}} \right) f_Y\left( y \right)dy}
& = & 2\pi x f_Y\left( x \right)\mathop {\lim }\limits_{m \to \infty } \mathop {\sup }\limits_{n \ge m} \int\limits {{{\left| {\frac{{b_n^\kappa K^{ft}\left( p \right)}}{{\overline {{g^{mt}}\left( {{p \mathord{\left/
 {\vphantom {t {{b_n}}}} \right.
 \kern-\nulldelimiterspace} {{b_n}}}} \right)} }}} \right|}^2}dp} \nonumber \\
&\le & \frac{{2\pi }}{C_1^2} x f_Y\left( x \right)\int{|p|^{2\kappa}\left| K^{ft}(p) \right|^2}dp \nonumber
\end{eqnarray}
at points of continuity $x > 0$ of $f_Y$, where
$${\phi _n}\left( y \right) = {\left| {\int\limits {{e^{ - ipy}}\frac{{b_n^\kappa K^{ft}\left( p \right)}}{{{g^{mt}}\left( {{p/{{b_n}}}} \right)}}dp} } \right|^2}$$ and
$$0 < {C_1} = \mathop {\lim }\limits_{u \to \infty } \mathop {\inf }\limits_{\left| p \right| \ge u} \left| {{g^{mt}}\left( p \right)} \right|{\left| p \right|^\kappa } \le \mathop {\lim }\limits_{u \to \infty } \mathop {\sup }\limits_{\left| p \right| \ge u} \left| {{g^{mt}}\left( p \right)} \right|{\left| p \right|^\kappa }  = {C_2}.$$
\end{lemma}	

\noindent \textit{Proof of Lemma \ref{lem1}}. 

We verify $\left( a \right)$ and use the same argument for $\left( b \right)$.\\
Let $h(\ln y) = y f_Y \left(y \right)$, we have
\begin{eqnarray}
	& & \mathop {\lim }\limits_{m \to \infty } \mathop {\inf }\limits_{n \ge m} \frac{1}{{{b_n}}}\int\limits {{\phi _n}\left( {\frac{{\ln x - \ln y}}{{{b_n}}}} \right) f_Y\left( y \right)dy}  - \frac{{2\pi }}{C_2^2} x f_Y\left( x \right)\int{|p|^{2\kappa}\left| K^{ft}(p) \right|^2}dp \nonumber \\
	&=& \mathop {\lim }\limits_{m \to \infty } \mathop {\inf }\limits_{n \ge m} \left\{ {\frac{1}{{{b_n}}}\int\limits {{\phi _n}\left( {\frac{{\ln x - \ln y}}{{{b_n}}}} \right)h\left( \ln y \right)d\left( \ln y \right)}  - \frac{{h\left(\ln x \right)}}{{{b_n}}}\int\limits {{\phi _n}\left( {\frac{{\ln x - \ln y}}{{{b_n}}}} \right)d\left( \ln y \right)} } \right. \nonumber \\
	& & \hspace{2cm} \left. { + \frac{{h\left(\ln x \right)}}{{{b_n}}}\int\limits {{\phi _n}\left( {\frac{{\ln x - \ln y}}{{{b_n}}}} \right)d\left( \ln y \right)} - \frac{{2\pi }}{C_2^2} {h\left(\ln x \right)}\int{|p|^{2\kappa}\left| K^{ft}(p) \right|^2}dp } \right\} \nonumber \\
	&=& \mathop {\lim }\limits_{m \to \infty } \mathop {\inf }\limits_{n \ge m} \left\{ {\int\limits {\left[ {h\left( {\ln x - \ln y} \right) - h\left(\ln x \right)} \right]\frac{1}{{{b_n}}}{\phi _n}\left( {\frac{\ln y}{{{b_n}}}} \right)d\left(\ln y\right)}} \right. \nonumber \\ 
	& & \hspace{2cm} \left. { + h\left(\ln x \right)\left[ {\int\limits {{\phi _n}\left(\ln y \right)d\left(\ln y\right) }  - \frac{{2\pi }}{C_2^2} \int{|p|^{2\kappa}\left| K^{ft}(p) \right|^2}dp } \right]} \right\}. \label{eq:61}
\end{eqnarray}
To reduce the limit of \eqref{eq:61}, we consider two steps:\\
\textbf{Step 1:} we prove $\displaystyle \int\limits {\left[ {h\left( {\ln x - \ln y} \right) - h\left(\ln x \right)} \right]\frac{1}{{{b_n}}}{\phi _n}\left( {\frac{\ln y}{{{b_n}}}} \right)d\left(\ln y\right) }  \to 0$ as $n \to \infty $.\\
We have 
$$\left| {\frac{{b_n^\kappa K^{ft}\left( p \right)}}{{{g^{mt}}\left( {{p/{b_n}}} \right)}}} \right| \le \frac{{b_n^\kappa \left| {K^{ft}\left( p \right)} \right|}}{{\mathop {\min }\limits_{\left| p \right| \le N} \left| {{g^{mt}}\left( p \right)} \right|}}{1_{\left\{ {\left| p \right| \le N{b_n}} \right\}}} + 2\frac{{{{\left| p \right|}^\kappa }\left| {K^{ft}\left( p \right)} \right|}}{{{C_1}}}{1_{\left\{ {\left| p \right| > N{b_n}} \right\}}},$$
where $N$ is large enough (but fixed) such that: ${\left| p \right|^\kappa }\left| {{g^{mt}}\left( p \right)} \right| \ge \frac{{{C_1}}}{2}$ when $\left| p \right| > N$. \\
We obtain
\begin{eqnarray}
	\left| {\int\limits {\frac{{b_n^\kappa K^{ft}\left( p \right)}}{{{g^{mt}}\left( {{p/{{b_n}}}} \right)}}dp} } \right| &\le & \int\limits {\left| {\frac{{b_n^\kappa K^{ft}\left( p \right)}}{{{g^{mt}}\left( {{p/{{b_n}}}} \right)}}} \right|dp } \nonumber \\
	&\le & \int\limits {\left( {\frac{{b_n^\kappa \left| {K^{ft}\left( p \right)} \right|}}{{\mathop {\min }\limits_{\left| p \right| \le N} \left| {{g^{mt}}\left( p \right)} \right|}}{1_{\left\{ {\left| p \right| \le N{b_n}} \right\}}} + 2\frac{{{{\left| p \right|}^\kappa }\left| {K^{ft}\left( p \right)} \right|}}{{{C_1}}}{1_{\left\{ {\left| p \right| > N{b_n}} \right\}}}} \right)dp} \nonumber \\
	&\le & \bigO\left( {b_n^\kappa } \right) + \frac{2}{{{C_1}}}\int\limits {{{\left| p \right|}^\kappa }\left| {K^{ft}\left( p \right)} \right|dp}  \le {D_1},\,\,\,\forall n \label{eq:62}
\end{eqnarray}
where ${D_1} > 0$ is a constant.
\begin{eqnarray}
	\int\limits {{e^{ - ip \ln y}}\frac{{b_n^\kappa K^{ft}\left( p \right)}}{{{g^{mt}}\left( {{p/{{b_n}}}} \right)}}dp}  
	&=&  - \frac{1}{{i \ln y}}\left( {\left. {{e^{ - ip \ln y}}\frac{{b_n^\kappa K^{ft}\left( p \right)}}{{{g^{mt}}\left( {{p/{{b_n}}}} \right)}}} \right|_{ - \infty }^{ + \infty } - \int\limits {{e^{ - ip \ln y}}{{\left( {\frac{{b_n^\kappa K^{ft}\left( p \right)}}{{{g^{mt}}\left( {{p/{{b_n}}}} \right)}}} \right)}^\prime }dp} } \right) \nonumber \\
	&=& - \frac{1}{{i \ln y}}{ \int\limits {{e^{ - ip \ln y}}\frac{{b_n^\kappa K^{ft}\left( p \right){{\left( {{g^{mt}}\left( {{p/{{b_n}}}} \right)} \right)}^\prime }}}{{{{\left[ {{g^{mt}}\left( {{p/{{b_n}}}} \right)} \right]}^2}}}dp} }. \nonumber
\end{eqnarray}
We deduce
$$ \left| {\int\limits {{e^{ - ip \ln y}}\frac{{b_n^\kappa K^{ft}\left( p \right)}}{{{g^{mt}}\left( {{p/{{b_n}}}} \right)}}dp} } \right| \le \frac{1}{{\left| \ln y \right|}} {\int\limits {\left| {\frac{{b_n^\kappa K^{ft}\left( p \right){{\left( {{g^{mt}}\left( {{p/{{b_n}}}} \right)} \right)}^\prime }}}{{{{\left[ {{g^{mt}}\left( {{p/{{b_n}}}} \right)} \right]}^2}}}} \right|dp} }. $$
By the similar argument as \eqref{eq:62}, we conclude
\begin{equation}
	\left| {\int\limits {{e^{ - ip \ln y}}\frac{{b_n^\kappa K^{ft}\left( p \right)}}{{{g^{mt}}\left( {{p/{{b_n}}}} \right)}}dp} } \right| \le \frac{{{D_2}}}{{\left| \ln y \right|}}\,\,\,\forall n\,\,\, \Rightarrow \mathop {\,\,\sup }\limits_n {\left| {\int\limits {{e^{ - ip \ln y}}\frac{{b_n^\kappa K^{ft}\left( p \right)}}{{{g^{mt}}\left( {{p/{{b_n}}}} \right)}}dp} } \right|^2} \le \frac{{D_2^2}}{{{{\left| \ln y \right|}^2}}}, \label{eq:63} 
\end{equation}
where ${D_2} > 0$ is a constant.

From \eqref{eq:62}, \eqref{eq:63} and using argument as in proof of Lemma 5.1 in Trong and Hung \cite{Trong2023}, we get
\begin{equation}
	\int\limits {\left[ {h\left( {\ln x - \ln y} \right) - h\left(\ln x \right)} \right]\frac{1}{{{b_n}}}{\phi _n}\left( {\frac{\ln y}{{{b_n}}}} \right)d\left(\ln y\right) }  \to 0 \text{ as } n \to \infty. \label{eq:64}
\end{equation}
\textbf{Step 2:} we prove $\displaystyle \mathop {\lim }\limits_{m \to \infty } \mathop {\inf }\limits_{n \ge m} \left\{ h\left(\ln x \right)\left[ {\int\limits {{\phi _n}\left(\ln y \right)d\left(\ln y\right) }  - \frac{{2\pi }}{C_2^2} \int{|p|^{2\kappa}\left| K^{ft}(p) \right|^2}dp } \right] \right\} \ge 0.$ We have
\begin{eqnarray}
	\int\limits {{\phi _n}\left(\ln y \right)d\left(\ln y\right) }
	&=& \int_{- \infty}^{+ \infty} {{\phi _n}\left( y \right)dy}  
	= \int_{- \infty}^{+ \infty} {{{\left| {\int\limits {{e^{ - ipy}}\frac{{b_n^\kappa K^{ft}\left( p \right)}}{{{g^{mt}}\left( {{p/{{b_n}}}} \right)}}dp} } \right|}^2}dy} \nonumber \\
	&=& \int_{- \infty}^{+ \infty} {{{\left| {{{\left[ {\frac{{b_n^\kappa K^{ft}\left( p \right)}}{{\overline {{g^{mt}}\left( {{p/{{b_n}}}} \right)} }}} \right]}^{ft}}\left( y \right)} \right|}^2}dy}  \nonumber \\
	&=& 2\pi \int_{-\infty}^{+\infty} {{{\left| {\frac{{b_n^\kappa K^{ft}\left( p \right)}}{{\overline {{g^{mt}}\left( {{p/{{b_n}}}} \right)} }}} \right|}^2}dp} \text{\,\, by Parseval’s identity.} \nonumber
\end{eqnarray}
We obtain
\begin{eqnarray}
	& & \mathop {\lim }\limits_{m \to \infty } \mathop {\inf }\limits_{n \ge m} \left\{ h\left(\ln x \right)\left[ {\int\limits {{\phi _n}\left(\ln y \right)d\left(\ln y\right) }  - \frac{{2\pi }}{C_2^2} \int{|p|^{2\kappa}\left| K^{ft}(p) \right|^2}dp } \right] \right\} \nonumber \\
	&=& h\left(\ln x \right)\left[ {\mathop {\lim }\limits_{m \to \infty } \mathop {\inf }\limits_{n \ge m} \int\limits {{\phi _n}\left(\ln y \right)d\left(\ln y\right)}  - \frac{{2\pi }}{C_2^2} \int{|p|^{2\kappa}\left| K^{ft}(p) \right|^2}dp } \right] \nonumber \\
	&=& h\left(\ln x \right)\left[ {\mathop {\lim }\limits_{m \to \infty } \mathop {\inf }\limits_{n \ge m} 2\pi \int\limits {{{\left| {\frac{{b_n^\kappa K^{ft}\left( p \right)}}{{\overline {{g^{mt}}\left( {{p/{{b_n}}}} \right)} }}} \right|}^2}dp}  - \frac{{2\pi }}{C_2^2} \int{|p|^{2\kappa}\left| K^{ft}(p) \right|^2}dp } \right] \nonumber \\
	&\ge & h\left(\ln x \right)\left[ {2\pi \int\limits {\mathop {\lim }\limits_{m \to \infty } \mathop {\inf }\limits_{n \ge m} \frac{{{{\left| p \right|}^{2\kappa }}\left| {K^{ft}\left( p \right)} \right|}}{{{{\left| {\overline {{g^{mt}}\left( {{p/{{b_n}}}} \right)} } \right|}^2}{{\left| {{p/{{b_n}}}} \right|}^{2\kappa }}}}dp}  - \frac{{2\pi }}{C_2^2} \int{|p|^{2\kappa}\left| K^{ft}(p) \right|^2}dp } \right] = 0. \quad \label{eq:65}
\end{eqnarray}
From \eqref{eq:61}, \eqref{eq:64} and \eqref{eq:65}, we conclude
\[\frac{{2\pi }}{C_2^2} x f_Y\left( x \right)\int{|p|^{2\kappa}\left| K^{ft}(p) \right|^2}dp  \le \mathop {\lim }\limits_{m \to \infty } \mathop {\inf }\limits_{n \ge m} \frac{1}{{{b_n}}}\int\limits {{\phi _n}\left( {\frac{{\ln x - \ln y}}{{{b_n}}}} \right) f_Y\left( y \right)dy}. \]
Similarly, we have
\begin{eqnarray}
	& & \mathop {\lim }\limits_{m \to \infty } \mathop {\inf }\limits_{n \ge m} \frac{1}{{{b_n}}}\int\limits {{\phi _n}\left( {\frac{{\ln x - \ln y}}{{{b_n}}}} \right)h\left(\ln y \right)d\left(\ln y \right)} \nonumber \\
	&=& \mathop {\lim }\limits_{m \to \infty } \mathop {\inf }\limits_{n \ge m} \left\{ {\frac{1}{{{b_n}}}\int\limits {{\phi _n}\left( {\frac{{\ln x - \ln y}}{{{b_n}}}} \right)h\left( \ln y \right)d\left( \ln y \right)}  - \frac{{h\left(\ln x \right)}}{{{b_n}}}\int\limits {{\phi _n}\left( {\frac{{\ln x - \ln y}}{{{b_n}}}} \right)d\left( \ln y \right)} } \right. \nonumber \\
	& & \hspace{8.6cm} \left. { + \frac{{h\left(\ln x \right)}}{{{b_n}}}\int\limits {{\phi _n}\left( {\frac{{\ln x - \ln y}}{{{b_n}}}} \right)d\left( \ln y \right)} } \right\} \nonumber \\
	&=& \mathop {\lim }\limits_{m \to \infty } \mathop {\inf }\limits_{n \ge m} \left\{ {\int\limits {\left[ {h\left( {\ln x - \ln y} \right) - h\left(\ln x \right)} \right]\frac{1}{{{b_n}}}{\phi _n}\left( {\frac{\ln y}{{{b_n}}}} \right)d\left(\ln y\right) }  + h\left(\ln x \right)\int\limits {{\phi _n}\left(\ln y \right)d\left(\ln y \right)} } \right\} \nonumber \\
	&=& 2\pi x f_Y\left( x \right)\mathop {\lim }\limits_{m \to \infty } \mathop {\inf }\limits_{n \ge m} \int\limits {{{\left| {\frac{{b_n^\kappa K^{ft}\left( p \right)}}{{\overline {{g^{mt}}\left( {{p/{{b_n}}}} \right)} }}} \right|}^2}dp}.  \,\,\, \nonumber 
\end{eqnarray}
\hfill ${\ensuremath{\square}}$

\noindent\textit{Proof of Theorem \ref{th:variance}}. First, we define
$${W_{{b_n}}}\left( x \right) := \frac{1}{{2\pi }}\int {{e^{ - ipx}}\frac{{K^{ft}\left( p \right)}}{{{g^{mt}}\left( {{p/{{b_n}}}} \right)}}dp}.$$
The estimator $\hat f_n$ can be represented as
\begin{eqnarray}
{\widehat f_n}\left( x \right) 
= \frac{1}{nx}\sum\limits_{j = 1}^n {\frac{1}{{{b_n}}}{W_{{b_n}}}\left( {\frac{{\ln x - \ln{Y_j}}}{{{b_n}}}} \right)}. \nonumber
\end{eqnarray}
Denote $${\mu _n} = \frac{1}{{{b_n}}}\mathbb{E}\left[ {{W_{{b_n}}}\left( {\frac{{\ln x - \ln{Y_j}}}{{{b_n}}}} \right)} \right] \text{ and } {U_{n,j}} = \frac{1}{{{b_n}}}{W_{{b_n}}}\left( {\frac{{\ln x - \ln{Y_j}}}{{{b_n}}}} \right) - {\mu_n}.$$
We have
\begin{eqnarray}
{\mu _n}
&=& \frac{1}{{{b_n}}}\mathbb{E}\left[ {\frac{1}{{2\pi }}\int\limits {{e^{ - ip\left( {\frac{{\ln x - \ln{Y_j}}}{{{b_n}}}} \right)}}\frac{{K^{ft}\left( p \right)}}{{{g^{mt}}\left( {{p/{{b_n}}}} \right)}}dp} } \right] = x \mathbb{E}\left[ {{{\widehat f}_n}\left( x \right)} \right] \to x{f_X}\left( x \right) \text{ as } n \to \infty, \nonumber 
\end{eqnarray}
where the last step follows by Theorem \ref{th:bias}. Then
\begin{eqnarray}
b_n^{1 + 2\kappa }\Var \left[ {{{U}_{n,j}}} \right] = b_n^{1 + 2\kappa }\mathbb{E}{\left| {{{U}_{n,j}}} \right|^2} &=& b_n^{1 + 2\kappa } \mathbb{E}{\left| {\frac{1}{{{b_n}}}\frac{1}{{2\pi }}\int {{e^{ - ip\left( {\frac{{\ln x - \ln{Y_j}}}{{{b_n}}}} \right)}}\frac{{K^{ft}\left( p \right)}}{{{g^{mt}}\left( {{p/{{b_n}}}} \right)}}dp}  - {\mu _n}} \right|^2} \nonumber \\
&=& \frac{1}{{{b_n}}}\mathbb{E}{\left| {\frac{1}{{2\pi }}\int\limits {{e^{ - ip\left( {\frac{{\ln x - \ln{Y_j}}}{{{b_n}}}} \right)}}\frac{{b_n^\kappa K^{ft}\left( p \right)}}{{{g^{mt}}\left( {{p/{{b_n}}}} \right)}}dp} } \right|^2} + \smallO\left( 1 \right) \nonumber \\
&=& \frac{1}{{4{\pi ^2}}}\frac{1}{{{b_n}}}\int\limits {{{\left| {\int\limits {{e^{ - ip\left( {\frac{{\ln x - \ln y}}{{{b_n}}}} \right)}}\frac{{b_n^\kappa K^{ft}\left( p \right)}}{{{g^{mt}}\left( {{p/{{b_n}}}} \right)}}dp} } \right|}^2} f_Y \left( y \right)dy}  + \smallO\left( 1 \right). \nonumber
\end{eqnarray}
By Lemma \ref{lem1}, we have
\begin{eqnarray}
x^2 \sigma _1^2\left( x \right) = \mathop {\lim }\limits_{m \to \infty } \mathop {\inf }\limits_{n \ge m} b_n^{1 + 2\kappa }\Var \left[ {{{U}_{n,j}}} \right] \le \mathop {\lim }\limits_{m \to \infty } \mathop {\sup }\limits_{n \ge m} b_n^{1 + 2\kappa }\Var \left[ {{{U}_{n,j}}} \right] = x^2 \sigma _2^2\left( x \right). \label{eq:66}
\end{eqnarray}
We have
\begin{eqnarray}
b_n^{1 + 2\kappa }\sum\limits_{l = 2}^n {\left| {\cov \left\{ {{{U}_{n,1}},{{U}_{n,l}}} \right\}} \right|}  
&=& b_n^{1 + 2\kappa }\sum\limits_{l = 2}^{{c_n}} {\left| {\cov \left\{ {{{U}_{n,1}},{{U}_{n,l}}} \right\}} \right|}  + b_n^{1 + 2\kappa }\sum\limits_{l = {c_n} + 1}^n {\left| {\cov \left\{ {{{U}_{n,1}},{{U}_{n,l}}} \right\}} \right|}  \nonumber \\
&=& {S_1} + {S_2}, \nonumber
\end{eqnarray} 
where ${c_n} \to \infty $ such that ${c_n}{b_n} \to 0$ as $n \to \infty $.

\noindent Let $Z_j = \ln Y_j$, $f_{Z_j} \left(u\right)$ and  $f_{Z_j, Z_k} \left(u, v\right) \left(1 \le j,k \le n\right)$ be the probability density function of $Z_j$ and the joint probability density function of $\left(Z_j, Z_k \right)$, respectively.  

\noindent For $2 \le l \le {c_n}$, we have in view of \eqref{eq:62} and \eqref{eq:63}
\begin{eqnarray}
& & \cov \left\{ {{{U}_{n,1}},{{U}_{n,l}}} \right\} 
= \frac{1}{{b_n^2}}{\cov} \left\{ {{W_{{b_n}}}\left( {\frac{{\ln x - \ln{Y_1}}}{{{b_n}}}} \right),{W_{{b_n}}}\left( {\frac{{\ln x - \ln{Y_l}}}{{{b_n}}}} \right)} \right\} \nonumber \\
&=& \frac{1}{{b_n^2}} \iint{{W_{{b_n}}}\left( {\frac{{\ln x - u}}{{{b_n}}}} \right){W_{{b_n}}}\left( {\frac{{\ln x - v}}{{{b_n}}}} \right)\left[ {{f_{{Z_1},{Z_l}}}\left( {u,v} \right) - f_{Z_1}\left( u \right)f_{Z_l}\left( v \right)} \right]dudv} \nonumber \\
&\le & \frac{{\mathop {\sup }\limits_{u,v,l} \left| {{f_{{Z_1},{Z_l}}}\left( {u,v} \right) - f_{Z_1}\left( u \right)f_{Z_l}\left( v \right)} \right|}}{{b_n^2}}{\left[ {\int\limits_{}^{} {\left| {{W_{{b_n}}}\left( {\frac{{\ln x - u}}{{{b_n}}}} \right)} \right|du} } \right]^2} \nonumber \\
& = & {\bigO\left( 1 \right)}{\left[ {\int\limits {\left| {{W_{{b_n}}}\left( u \right)} \right|du} } \right]^2} \le \frac{{\bigO\left( 1 \right)}}{{{b^{2\kappa }}}}. \nonumber
\end{eqnarray}
Then
$$ {S_1} = b_n^{1 + 2\kappa }\sum\limits_{l = 2}^{{c_n}} {\left| {\cov \left\{ {{{U}_{n,1}},{{U}_{n,l}}} \right\}} \right|}  \le b_n^{1 + 2\kappa }{c_n}\frac{{\bigO\left( 1 \right)}}{{{b^{2\kappa }}}} = \bigO\left( {{b_n}{c_n}} \right) \to 0 \,\,{\mathop{\rm as}\nolimits} \,\,n \to \infty. $$
For ${c_n} + 1 \le l \le n$, note that the process $\left( {{Y_j}} \right)$ is also strongly mixing (see Bradley \cite{Bradley2005}) with strong mixing coefficient
${\alpha _Y}\left( k \right) \le \bigO\left( {{1 \mathord{\left/
 {\vphantom {1 {{k^{2 + \delta }}}}} \right.
 \kern-\nulldelimiterspace} {{k^{2 + \delta }}}}} \right)$. This condition is equivalent to for some $\nu  > 2$ and $a > 1 - {2 \mathord{\left/
 {\vphantom {2 \nu }} \right.
 \kern-\nulldelimiterspace} \nu },$
\begin{eqnarray}
\sum\limits_{k = 1}^\infty  {{k^a}{{\left[ {{\alpha _Y}\left( k \right)} \right]}^{1 - {2 \mathord{\left/
 {\vphantom {2 \nu }} \right.
 \kern-\nulldelimiterspace} \nu }}}}  < \infty. \label{eq:67}
\end{eqnarray} 
 
\noindent Applying Davydov's Lemma (see Hall and Heyde \cite{Hall1980}, Corollary A.2), we find
\[\left| {\cov \left\{ {{{U}_{n,1}},{{U}_{n,l}}} \right\}} \right| \le \frac{8}{{b_n^2}}{\left[ {{\alpha _Y}\left( {l - 1} \right)} \right]^{1 - {2 \mathord{\left/
 {\vphantom {2 \nu }} \right.
 \kern-\nulldelimiterspace} \nu }}}{\left\{ {\mathbb{E}\left[ {{{\left| {{W_{{b_n}}}\left( {\frac{{\ln x - \ln{Y_1}}}{{{b_n}}}} \right)} \right|}^\nu }} \right]} \right\}^{{2 \mathord{\left/
 {\vphantom {2 \nu }} \right.
 \kern-\nulldelimiterspace} \nu }}}.\]
Note that
$$\frac{1}{{{b_n}}}{\mathbb{E}\left[ {{{\left| {{W_{{b_n}}}\left( {\frac{{\ln x - \ln{Y_1}}}{{{b_n}}}} \right)} \right|}^\nu }} \right]} \le \bigO\left( 1 \right)\int\limits {{{\left| {{W_{{b_n}}}\left( u \right)} \right|}^\nu }du}  \le \frac{{O\left( 1 \right)}}{{b_n^{\nu \kappa }}},$$
where the last step follows by $\eqref{eq:62}$ and $\eqref{eq:63}$.\\
Hence,
$$ \left| {\cov \left\{ {{{U}_{n,1}},{{U}_{n,l}}} \right\}} \right| \le \frac{{\bigO\left( 1 \right)}}{{b_n^{2 + 2\kappa  - {2 \mathord{\left/
 {\vphantom {2 \nu }} \right.
 \kern-\nulldelimiterspace} \nu }}}}{\left[ {\alpha_{Y} \left( {l - 1} \right)} \right]^{1 - {2 \mathord{\left/
 {\vphantom {2 \nu }} \right.
 \kern-\nulldelimiterspace} \nu }}}$$
and
\begin{eqnarray}
{S_2} 
&=& b_n^{1 + 2\kappa }\sum\limits_{l = {c_n} + 1}^n {\left| {\cov \left\{ {{{U}_{n,1}},{{U}_{n,l}}} \right\}} \right|}  \le \frac{{\bigO\left( 1 \right)}}{{b_n^{1 - {2 \mathord{\left/
 {\vphantom {2 \nu }} \right.
 \kern-\nulldelimiterspace} \nu }}}}\sum\limits_{l = {c_n} + 1}^n {{{\left[ {\alpha_{Y} \left( {l - 1} \right)} \right]}^{1 - {2 \mathord{\left/
 {\vphantom {2 \nu }} \right.
 \kern-\nulldelimiterspace} \nu }}}} \nonumber \\
&\le & \frac{{\bigO\left( 1 \right)}}{{c_n^ab_n^{1 - {2 \mathord{\left/
 {\vphantom {2 \nu }} \right.
 \kern-\nulldelimiterspace} \nu }}}}\sum\limits_{j = {c_n}}^\infty  {{j^a}{{\left[ {\alpha_{Y} \left( j \right)} \right]}^{1 - {2 \mathord{\left/
 {\vphantom {2 \nu }} \right.
 \kern-\nulldelimiterspace} \nu }}}}. \nonumber
\end{eqnarray}
Choose ${c_n} = {\left( {\dfrac{1}{{b_n^{1 - {2/\nu }}}}} \right)^{1/a}}$ then for $a > 1 - {2/\nu }$ we obtain ${c_n}{b_n} \to 0$ as required. Then
\[
{S_2} \le \bigO\left( 1 \right)\sum\limits_{j = {c_n}}^\infty  {{j^a}{{\left[ {\alpha_{Y} \left( j \right)} \right]}^{1 - {2/\nu }}}}  \to 0 \text{\,\,as\,\,}n \to \infty \text{ follows by \eqref{eq:67}. } 
\]
Therefore, 
\begin{equation}
\mathop {\lim }\limits_{n \to \infty } b_n^{1 + 2\kappa }\sum\limits_{l = 2}^n {\left| {\cov \left\{ {{{U}_{n,1}},{{U}_{n,l}}} \right\}} \right|}  = 0. \label{eq:68}
\end{equation}
We have
\begin{eqnarray}
& & n\Var \left[ {{{\widehat f}_n}\left( x \right)} \right] = n\Var \left\{ {\frac{1}{nx}\sum\limits_{j = 1}^n {\frac{1}{{{b_n}}}{W_{{b_n}}}\left( {\frac{{\ln x - \ln{Y_j}}}{{{b_n}}}} \right)} } \right\} \nonumber \\
&=& \frac{1}{x^2} \Var \left[ {\frac{1}{{{b_n}}}{W_{{b_n}}}\left( {\frac{{\ln x - \ln{Y_j}}}{{{b_n}}}} \right)} \right] \nonumber \\
& & + \frac{2}{nx^2}\sum\limits_{j = 2}^n {\left[ {n - (j-1)} \right]\cov \left\{ {\frac{1}{{{b_n}}}{W_{{b_n}}}\left( {\frac{{\ln x - \ln{Y_1}}}{{{b_n}}}} \right),\frac{1}{{{b_n}}}{W_{{b_n}}}\left( {\frac{{\ln x - \ln{Y_j}}}{{{b_n}}}} \right)} \right\}} \nonumber \\
&=& \frac{1}{x^2} \Var \left[ {\frac{1}{{{b_n}}}{W_{{b_n}}}\left( {\frac{{\ln x - \ln{Y_j}}}{{{b_n}}}} \right)} \right] \nonumber \\
& & + \frac{2}{x^2} \sum\limits_{j = 2}^n {\left( {1 - \frac{j-1}{n}} \right){\mathop{\rm cov}} \left\{ {\frac{1}{{{b_n}}}{W_{{b_n}}}\left( {\frac{{\ln x - \ln{Y_1}}}{{{b_n}}}} \right),\frac{1}{{{b_n}}}{W_{{b_n}}}\left( {\frac{{\ln x - \ln{Y_j}}}{{{b_n}}}} \right)} \right\}} \nonumber \\
&=& \frac{1}{x^2} \Var \left[ {{{U}_{n,j}}} \right] + \frac{2}{x^2} \sum\limits_{j = 2}^n {\left( {1 - \frac{j-1}{n}} \right)\cov \left\{ {{{U}_{n,1}},{{U}_{n,j}}} \right\}}. \nonumber
\end{eqnarray}
By $\eqref{eq:66}$ and $\eqref{eq:68}$, we obtain
\[\sigma _1^2\left( x \right) = \mathop {\lim }\limits_{m \to \infty } \mathop {\inf }\limits_{n \ge m} nb_n^{1 + 2\kappa }\Var \left[ {{{\widehat f}_n}\left( x \right)} \right] \le \mathop {\lim }\limits_{m \to \infty } \mathop {\sup }\limits_{n \ge m} nb_n^{1 + 2\kappa }\Var \left[ {{{\widehat f}_n}\left( x \right)} \right] = \sigma _2^2\left( x \right).\]
This completes the proof of Theorem \ref{th:variance}. \hfill {\ensuremath{\square}}

\subsection{Proof of Theorem \ref{th: converge-a.s-estimators}}

To prove Theorem \ref{th: converge-a.s-estimators}, we first prove the following lemma.

% LEMMA 2
% LEMMA 2

\begin{lemma}\label{lem2}
Let $\left( {{Y_j}} \right)$ be a stationary strongly mixing process on the probability space $\left( {\Omega ,\mathcal{A},\pr} \right)$ with the mixing coefficient  satisfying ${\alpha}\left( k \right) \le a{\rho ^k} \left(a > 0,\,0 < \rho  < 1 \right)$. Let ${Z_j} = {\psi _1}\left( {{Y_j}} \right) + i.{\psi _2}\left( {{Y_j}} \right)$ where $\psi_1$ and $\psi_2$ are real-valued Borel measurable function. Assume that $\mathbb{E}\left( {{Z_j}} \right) = 0$ and $\left| {{Z_j}} \right| \le d$. Let an  integer $n \ge 1$ be given. Then  for all $\varepsilon  > 0$, 
\[\mathbb{P}\left\{ {\frac{1}{n}\left| {\sum\limits_{j = 1}^n {{Z_j}} } \right| \ge \varepsilon } \right\} \le 2\left( {1 + 4a{e^{ - 2}}L} \right)\exp \left( { - \frac{{{\varepsilon ^2}n}}{{4\left( {\mathbb{E}{{\left| {{Z_1}} \right|}^2} + {{\sqrt 2 \varepsilon d} \mathord{\left/
 {\vphantom {{\sqrt 2 \varepsilon d} 6}} \right.
 \kern-\nulldelimiterspace} 6}} \right)}}} \right),\]
 where $L \ge 1$ is a constant.
\end{lemma}

\noindent \textit{Proof of Lemma \ref{lem2}}. We have
\begin{eqnarray}
\left\{ {\frac{1}{n}\left| {\sum\limits_{j = 1}^n {{Z_j}} } \right| \ge \varepsilon } \right\} &=& \left\{ {\frac{1}{n}\sqrt {{{\left[ {\sum\limits_{j = 1}^n {{\psi _1}\left( {{Y_j}} \right)} } \right]}^2} + {{\left[ {\sum\limits_{j = 1}^n {{\psi _2}\left( {{Y_j}} \right)} } \right]}^2}}  \ge \varepsilon } \right\} \nonumber \\
&\subset& \left\{ {\frac{1}{n}\left| {\sum\limits_{j = 1}^n {{\psi _1}\left( {{Y_j}} \right)} } \right| \ge \frac{{\sqrt 2 }}{2}\varepsilon } \right\} \cup \left\{ {\frac{1}{n}\left| {\sum\limits_{j = 1}^n {{\psi _2}\left( {{Y_j}} \right)} } \right| \ge \frac{{\sqrt 2 }}{2}\varepsilon } \right\}. \nonumber
\end{eqnarray}
We deduce
\[ \mathbb{P}\left\{ {\frac{1}{n}\left| {\sum\limits_{j = 1}^n {{Z_j}} } \right| \ge \varepsilon } \right\} \le \mathbb{P}\left\{ {\frac{1}{n}\left| {\sum\limits_{j = 1}^n {{\psi _1}\left( {{Y_j}} \right)} } \right| \ge \frac{{\sqrt 2 }}{2}\varepsilon } \right\} + \mathbb{P}\left\{ {\frac{1}{n}\left| {\sum\limits_{j = 1}^n {{\psi _2}\left( {{Y_j}} \right)} } \right| \ge \frac{{\sqrt 2 }}{2}\varepsilon } \right\}.\]
Note that\\
\hspace*{.5cm}$\mathbb{E}\left( {{Z_j}} \right) = \mathbb{E}\left[ {{\psi _1}\left( {{Y_j}} \right) + i.{\psi _2}\left( {{Y_j}} \right)} \right] = 0$. This implies $\mathbb{E}\left[ {{\psi _1}\left( {{Y_j}} \right)} \right] = \mathbb{E}\left[ {{\psi _2}\left( {{Y_j}} \right)} \right] = 0.$\\
\hspace*{.5cm}$\left| {{\psi _k}\left( {{Y_j}} \right)} \right| \le \left| {{Z_j}} \right| \le d,\,\,\,k = 1,2.$\\
From Theorem 2.1 of N'drin and Hili \cite{Ndrin2013}, we have 
\[\mathbb{P}\left\{ {\frac{1}{n}\left| {\sum\limits_{j = 1}^n {{\psi _k}\left( {{Y_j}} \right)} } \right| \ge \frac{\sqrt{2}}{2}\varepsilon } \right\} \le \left( {1 + 4a{e^{ - 2}}L} \right)\exp \left( { - \frac{{{\varepsilon ^2}n}}{{4\left( {\mathbb{E}{{\left| {{\psi_k \left( Y_1 \right)}} \right|}^2} + {{\sqrt{2}\varepsilon d}/6}} \right)}}} \right).\]
Therefore,
\[
\mathbb{P}\left\{ {\frac{1}{n}\left| {\sum\limits_{j = 1}^n {{Z_j}} } \right| \ge \varepsilon } \right\} \le 2\left( {1 + 4a{e^{ - 2}}L} \right)\exp \left( { - \frac{{{\varepsilon ^2}n}}{{4\left( {\mathbb{E}{{\left| {{Z_1}} \right|}^2} + {{\sqrt 2 \varepsilon d}/6}} \right)}}} \right). 
\]
\hfill {\ensuremath{\square}}

\noindent \textit{Proof of Theorem \ref{th: converge-a.s-estimators}}. Notice that
\begin{eqnarray}
\left| {{{\widehat f}_n}\left( x \right) - {f_X }\left( x \right)} \right| 
&=& \left| {{{\widehat f}_n}\left( x \right) - \mathbb{E}{{\widehat f}_n}\left( x \right) + \mathbb{E}\widehat {{f_n}}\left( x \right) - {f_X }\left( x \right)} \right| \nonumber \\
&\le & \left| {{{\widehat f}_n}\left( x \right) - \mathbb{E}{{\widehat f}_n}\left( x \right)} \right| + \left| {\mathbb{E}{{\widehat f}_n}\left( x \right) - {f_X }\left( x \right)} \right|. \nonumber
\end{eqnarray}

Then we get
\begin{eqnarray}
\left| {{{\widehat f}_n}\left( x \right) - \mathbb{E}{{\widehat f}_n}\left( x \right)} \right| 
&=& \left| {\frac{1}{{2\pi }}\int\limits {{x^{ - 1 - ip}}\frac{K^{ft}\left( {p{b_n}} \right)}{{{g^{mt}}\left( p \right)}}\frac{1}{n}\sum\limits_{j = 1}^n {\left( {{Y_j^{ip}} - \mathbb{E}{Y_j^{ip}}} \right)} dp} } \right| \nonumber \\
&=& \frac{1}{{nb_n^{1 + \kappa }}}\left| {\sum\limits_{j = 1}^n {\frac{1}{{2\pi }}\int\limits {{x^{ - 1 - ip}}\frac{{b_n^{1 + \kappa }K^{ft}\left( {p{b_n}} \right)}}{{{g^{mt}}\left( p \right)}}{\left( {{Y_j^{ip}} - \mathbb{E}{Y_j^{ip}}} \right)}dp} } } \right| \nonumber \\
&=& \frac{1}{{nb_n^{1 + \kappa }}}\left| {\sum\limits_{j = 1}^n {{Z_j}} } \right|, \nonumber
\end{eqnarray}
where $\displaystyle {Z_j} = {\frac{1}{{2\pi }}\int\limits {{x^{ - 1 - ip}}\frac{{b_n^{1 + \kappa }K^{ft}\left( {p{b_n}} \right)}}{{{g^{mt}}\left( p \right)}}{\left( {{Y_j^{ip}} - \mathbb{E}{Y_j^{ip}}} \right)}dp} }$. 

We verify that $Z_j$ satisfy the assumptions of Lemma \ref{lem2}. We first have ${\left( {{Z_j}} \right)_{j = \overline {1,n} }}$ are identical distributions, $\mathbb{E}\left[ {{Z_j}} \right] = 0,\forall j$. For $\forall x > 0$, we have
\begin{eqnarray*}
\left| {{Z_j}} \right| 
&=& \left| {\frac{1}{{2\pi }}\int\limits {{x^{ - 1 - ip}}\frac{{b_n^{1 + \kappa }K^{ft}\left( {p{b_n}} \right)}}{{{g^{mt}}\left( p \right)}}{\left( {{Y_j^{ip}} - \mathbb{E}{Y_j^{ip}}} \right)}dp} } \right| \nonumber \\
&\le & \frac{1}{{2\pi x }}\int\limits {\left| {{x^{ - ip}}} \right|\left| {\frac{{b_n^{1 + \kappa }K^{ft}\left( {p{b_n}} \right)}}{{{g^{mt}}\left( p \right)}}} \right|\left| {{Y_j^{ip}} - \mathbb{E}{Y_j^{ip}}} \right|dp} \nonumber \\
&\le & \frac{1}{\pi x }\int\limits {\left| {\frac{{b_n^\kappa K^{ft}\left( p \right)}}{{{g^{mt}}\left( {{p/{{b_n}}}} \right)}}} \right|dp}   \nonumber \\
&\le & \frac{1}{{c \pi x}}\int\limits {{{\left| p \right|}^\kappa } {K^{ft}\left( p \right)} dp } = B_1, \nonumber 
\end{eqnarray*}
where ${B_1} > 0$ is a constant. So $Z_j$ satisfy assumptions of Lemma \ref{lem2}. For all $\varepsilon  > 0$, we get
\begin{eqnarray}
\mathbb{P}\left\{ {\frac{1}{{nb_n^{1 + \kappa }}}\left| {\sum\limits_{j = 1}^n {{Z_j}} } \right| \ge \varepsilon } \right\} 
&=& \mathbb{P}\left\{ {\frac{1}{n}\left| {\sum\limits_{j = 1}^n {{Z_j}} } \right| \ge \varepsilon b_n^{1 + \kappa }} \right\} \nonumber \\
&\le & 2\left( {1 + 8a{e^{ - 2}}L} \right)\exp \left( { - \frac{{{\varepsilon ^2}nb_n^{2 + 2\kappa }}} {{4\left[ {\mathbb{E}{{\left| {{Z_1}} \right|}^2} + {{\sqrt 2 \varepsilon b_n^{1 + \kappa } }}} B_1/6 \right] }}} \right) \nonumber
\end{eqnarray}
We have
\begin{eqnarray*}
\mathbb{E}{\left| {{Z_1}} \right|^2}
&=& \mathbb{E}{\left| {\frac{1}{{2\pi }}\int\limits {{x^{ - 1 - ip}}\frac{{b_n^{1 + \kappa }K^{ft}\left( {p{b_n}} \right)}}{{{g^{mt}}\left( p \right)}}\left( {{{Y_1}^{ip}} - \mathbb{E}{{Y_1}^{ip}}} \right)dp} } \right|^2} \nonumber \\
&\le & \frac{{b_n^{2 + 2\kappa }}}{{4{\pi ^2} x^2 }}\mathbb{E}{\left| {\int\limits {{e^{ip\left( {\ln{Y_1} - \ln x} \right)}}\frac{{K^{ft}\left( {p{b_n}} \right)}}{{{g^{mt}}\left( p \right)}}dp} } \right|^2} \nonumber \\
&=& \frac{b_n}{{4{\pi ^2} x^2 }} \frac{1}{b_n} \int {{{\left| {\int\limits {{e^{-ip\left( \frac{{\ln x - \ln y}}{b_n} \right)}}\frac{{b_n^{ \kappa }}{K^{ft}\left( {p} \right)}}{{{g^{mt}}\left( p/b_n \right)}}dp} } \right|}^2} f_Y\left( y \right)dy} \nonumber \\
&\le & \frac{{b_n}}{2\pi x C_1^2} f_Y\left( x \right)\int{|p|^{2\kappa}\left| K^{ft}(p) \right|^2}dp = b_n B_2, \nonumber 
\end{eqnarray*}
where ${B_2} > 0$ is a constant and the last step follows by Lemma \ref{lem1}. Moreover, the process $\left( {{Y_j}} \right)$ is also strongly mixing (see \cite{Bradley2005}) with strong mixing coefficient
\[{\alpha _Y}\left( k \right) \le {\alpha _X}\left( k \right) + {\alpha _U }\left( k \right) \le 2a{ {\rho}^k}.\]
This implies
\begin{eqnarray}
\mathbb{P}\left\{ {\frac{1}{{nb_n^{1 + \kappa }}}\left| {\sum\limits_{j = 1}^n {{Z_j}} } \right| \ge \varepsilon } \right\} 
&\le & 2\left( {1 + 8a{e^{ - 2}}L} \right)\exp \left( { - \frac{{{\varepsilon ^2}nb_n^{1 + 2\kappa }}}{{4\left[ {{B_2} + {{\sqrt 2 \varepsilon b_n^{ \kappa } }}} B_1/6 \right]}}} \right) \to 0 \nonumber
\end{eqnarray}
as $nb_n^{1 + 2\kappa } \to \infty.$

Therefore, by Borel-Cantelli's lemma, we deduce that ${\widehat f_n}\left( x \right) - \mathbb{E}{\widehat f_n}\left( x \right) \to 0$ a.s.  Thus, ${\widehat f_n}\left( x \right) \to {f_X}\left( x \right)$ a.s. by Theorem \ref{th:bias}. \hfill {\ensuremath{\square}}

\subsection{Proof of Theorem \ref{th:asymptotic-normality-estimators}}

We will prove \eqref{eq:33} and use the same argument for \eqref{eq:34}. Let ${V_{n,j}} = {b^{1/2 + \kappa }}{U_{n,j}}$, we have
\begin{eqnarray}
{n^{1/2}}{b^{{1/2} + \kappa }}\left( {{{\widehat f}_n}\left( x \right) - \mathbb{E}{{\widehat f}_n}\left( x \right)} \right) 
&=& {n^{1/2}}{b^{{1/2} + \kappa }}\frac{1}{nx}\sum\limits_{j = 1}^n {\left( {\frac{1}{{{b_n}}}{W_{{b_n}}}\left( {\frac{{\ln x - \ln{Y_j}}}{{{b_n}}}} \right) - {\mu _n}} \right)} \nonumber \\
&=& \frac{1}{{\sqrt n }}\sum\limits_{j = 1}^n {{b^{{1/2} + \kappa }}{{U}_{n,j}}}  = \frac{1}{{\sqrt n }}\sum\limits_{j = 1}^n {{V_{n,j}}}  = \,\frac{1}{{\sqrt n }}{S_n}. \nonumber
\end{eqnarray}
Note that the process $\left( {{Y_j}} \right)$ is also strongly mixing (see \cite{Bradley2005}) with strong mixing coefficient
${\alpha _Y}\left( k \right) \le \bigO\left( {{1/{{k^{2 + \delta }}}}} \right)$. And due to ${\left( {{n/{{b_n}}}} \right)^{1/2}}\alpha_U \left( {{s_n}} \right) \to 0\,\,{\rm{as}}\,\,n \to \infty $, we also have ${\left( {{n/{{b_n}}}} \right)^{1/2}}\alpha_Y \left( {{s_n}} \right) \to 0\,\,{\rm{as}}\,\,n \to \infty $. The conditions ${s_n} = \smallO\left( {{{\left( {n{b_n}} \right)}^{{1/2}}}} \right)$ and ${\left( {{n/{{b_n}}}} \right)^{1/2}}\alpha_Y \left( {{s_n}} \right) \to 0\,\,{\rm{as}}\,\,n \to \infty $ imply that there exists a sequence of integers $\left\{ {{q_n}} \right\},\,{q_n} \to \infty $, such that
 $${q_n}{s_n} = \smallO \left( {{{\left( {n{b_n}} \right)}^{{1/2}}}} \right)\,\,{\rm{and}}\,\,{q_n}{\left( {n \mathord{\left/
 {\vphantom {n {{b_n}}}} \right.
 \kern-\nulldelimiterspace} {{b_n}}} \right)^{1/2}}\alpha_Y \left( {{s_n}} \right) \to 0\,\,{\rm{as}}\,\,n \to \infty .$$
Let ${r_n} = \left\lfloor {\dfrac{{{{\left( {n{b_n}} \right)}^{1/2}}}}{{{q_n}}}} \right\rfloor $, we have some following properties as $n \to \infty $
\begin{eqnarray}
\frac{{{s_n}}}{{{r_n}}} \to 0 \quad \left( a \right), \quad \quad \frac{{{r_n}}}{n} \to 0 \quad \left( b \right), \quad \quad \frac{{{r_n}}}{{{{\left( {n{b_n}} \right)}^{{1 \mathord{\left/
 {\vphantom {1 2}} \right.
 \kern-\nulldelimiterspace} 2}}}}} \to 0 \quad \left( c \right), \quad \quad \frac{n}{{{r_n}}}\alpha_Y \left( {{s_n}} \right) \to 0 \quad \left( d \right) \label{eq:69}
\end{eqnarray}
Now, we partition the set $\left\{ {1,2,...,n} \right\}$ into $2k + 1$ subsets: $\left\{j(r+s)+1,\ldots,j(r+s)+r\right\}$, $\left\{j(r+s)+r+1,\ldots,(j+1)(r+s)\right\}$ for $0 \le j \le k-1$ and $\left\{k(r+s)+1,\ldots,n\right\}$, where $r = {r_n}, s = {s_n}$ and 
\begin{eqnarray}
k = {k_n} = \left\lfloor {\frac{n}{{{r_n} + {s_n}}}} \right\rfloor \label{eq:610}
\end{eqnarray}
For $0 \le j \le k - 1$, define the random variables 
\begin{eqnarray}
{\eta _j} = \sum\limits_{i = j\left( {r + s} \right) + 1}^{j\left( {r + s} \right) + r} {{V_{n,i}}} \quad \left( a \right), \quad {\xi _j} = \sum\limits_{i = j\left( {r + s} \right) + r + 1}^{\left( {j + 1} \right)\left( {r + s} \right)} {{V_{n,i}}} \quad \left( b \right), \quad {\zeta _k} = \sum\limits_{i = k\left( {r + s} \right) + 1}^n {{V_{n,i}}} \quad \left( c \right)  \label{eq:611}
\end{eqnarray}
Write ${S_n} = \sum\limits_{j = 0}^{k - 1} {{\eta _j}}  + \sum\limits_{j = 0}^{k - 1} {{\xi _j}}  + {\zeta _k} = {S'_n} + {S''_n} + {S'''_n}$, we will prove that, as $n \to \infty $,
\begin{eqnarray}
\begin{array}{*{20}{cc}}
{\frac{1}{n}\mathbb{E}{{\left( {{{S''}_n}} \right)}^2} \to 0,\,\,\,\,\,\,\frac{1}{n}\mathbb{E}{{\left( {{{S'''}_n}} \right)}^2} \to 0}  \quad \left( a \right)
&\quad &{\left| {\mathbb{E}\left( {{e^{ip{{S'}_n}}}} \right) - \prod\limits_{j = 0}^{k - 1} {\mathbb{E}\left( {{e^{ip{\eta _j}}}} \right)} } \right| \to 0} \quad \left( b \right)\\
{\mathop {\inf }\limits_n \left[ {\frac{1}{n}\sum\limits_{j = 0}^{k - 1} {\mathbb{E}\left( {\eta _j^2} \right)} } \right] \to \sigma _1^2\left( x \right)\,\,\,\,\,\,\,\,\,\,\,\,\,\,\,} \quad \left( c \right)
&\quad &{\frac{1}{n}\sum\limits_{j = 0}^{k - 1} {\mathbb{E}\left[ {\eta _j^2{I_{\left\{ {\left| {{\eta _j}} \right| > \varepsilon {\sigma _1}\left( x \right)\sqrt n } \right\}}}} \right]}  \to 0} \quad \left( d \right) \label{eq:612}
\end{array}
\end{eqnarray}
Note that $\left(\ref{eq:612}a\right)$ implies that ${S''_n},\,{S'''_n}$ are asymptotically negligible, $\left(\ref{eq:612}b\right)$ show that the summands $\left\{ {{\eta _j}} \right\}$ in  ${S'_n}$ are asymptotically independent, $\left(\ref{eq:612}c\right)$ and $\left(\ref{eq:612}d\right)$ are the standard Lindeberg-Feller conditions for asymptotic normality of ${S'_n}$ under independence.\\
Since $\mathbb{E}\left( {{V_{n,i}}} \right) = 0$, we have $\mathbb{E}\left( {{\eta _j}} \right) = \mathbb{E}\left( {{\xi _j}} \right) = \mathbb{E}\left( {{\zeta _k}} \right) = 0$. Consider first
$$\mathbb{E}{\left( {{{S''}_n}} \right)^2} = \Var \left[ {\sum\limits_{j = 0}^{k - 1} {{\xi _j}} } \right] = \sum\limits_{j = 0}^{k - 1} {\Var \left[ {{\xi _j}} \right]}  + \sum\limits_{i = 0}^{k - 1} {\sum\limits_{\scriptstyle j = 0\hfill\atop
\scriptstyle j \ne i\hfill}^{k - 1} {\cov \left\{ {{\xi _i},{\xi _j}} \right\}} }  = {F_1} + {F_2}.$$
With ${m_j} = j\left( {r + s} \right) + r$, by $\left(\ref{eq:611}b\right)$ we have
\begin{eqnarray}
\Var \left[ {{\xi _j}} \right] 
&=& \sum\limits_{{i_1} = 1}^s {\sum\limits_{{i_2} = 1}^s {\cov \left\{ {{V_{n,{m_j} + {i_1}}},{V_{n,{m_j} + {i_2}}}} \right\}} }  = s\Var \left[ {{V_{n,1}}} \right] + 2s\sum\limits_{i = 2}^s {\left( {1 - \frac{i}{s}} \right)\cov \left\{ {{V_{n,1}},{V_{n,i}}} \right\}} \nonumber \\
&\le & s\sigma _2^2\left( x \right)\left\{ {1 + \smallO \left( 1 \right)} \right\}, \nonumber
\end{eqnarray}
the last step follows by Theorem \ref{th:variance} and we obtain $ {F_1} \le \bigO\left( 1 \right)ks$. \\
Next, consider the term  ${F_2} = \sum\limits_{i = 0}^{k - 1} {\sum\limits_{\scriptstyle j = 0\hfill\atop
\scriptstyle j \ne i\hfill}^{k - 1} {\cov \left\{ {{\xi _i},{\xi _j}} \right\}} }  = \sum\limits_{i = 0}^{k - 1} {\sum\limits_{\scriptstyle j = 0\hfill\atop
\scriptstyle j \ne i\hfill}^{k - 1} {\sum\limits_{{l_1} = 1}^s {\sum\limits_{{l_2} = 1}^s {\cov \left\{ {{V_{n,{m_i} + {l_1}}},{V_{n,{m_j} + {l_2}}}} \right\}} } } } $ \\
Since $i \ne j, \left| {{m_j} + {l_2} - {m_i} - {l_1}} \right| \ge r$, so that
$${F_2} \le 2\sum\limits_{{l_1} = 1}^{n - r} {\sum\limits_{{l_2} = r + {l_1}}^n {\left| {\cov \left\{ {{V_{n,{l_1}}},{V_{n,{l_2}}}} \right\}} \right|} }  \le 2n\sum\limits_{j = r + 1}^n {\left| {\cov \left\{ {{V_{n,1}},{V_{n,j}}} \right\}} \right|}  = \smallO \left( n \right).$$
Therefore, we obtain in view of $\left(\ref{eq:69}a\right)$ and \eqref{eq:610}
$$\frac{1}{n}\mathbb{E}{\left( {{{S''}_n}} \right)^2} = \frac{{{F_1} + {F_2}}}{n} \le \frac{{\bigO\left( 1 \right)ks + \smallO \left( n \right)}}{n} \to 0\,\,{\rm{as}}\,\,n \to \infty.$$
Using a similar argument for the latter inequality, we have in view of $\left(\ref{eq:69}a\right)$ and $\left(\ref{eq:69}b\right)$ that
\begin{eqnarray}
\frac{1}{n}\mathbb{E}{\left( {{{S'''}_n}} \right)^2} 
&=& \frac{1}{n}\sum\limits_{j = 1}^{n - k\left( {r + s} \right)} {\Var \left[ {{V_{n,k\left( {r + s} \right) + j}}} \right]}  + \frac{1}{n}\sum\limits_{i = 1}^{n - k\left( {r + s} \right)} {\sum\limits_{\scriptstyle j = 1\hfill\atop
\scriptstyle j \ne i\hfill}^{n - k\left( {r + s} \right)} {\cov \left\{ {{V_{n,k\left( {r + s} \right) + i}},{V_{n,k\left( {r + s} \right) + j}}} \right\}} } \nonumber \\
&\le & \frac{C}{n}\left[ {n - k\left( {r + s} \right)} \right] + \frac{{\smallO \left( n \right)}}{n} \le \frac{{C\left( {r + s} \right)}}{n} + \frac{{\smallO \left( n \right)}}{n} \to 0\,\,{\rm{as}}\,\,n \to \infty. \nonumber
\end{eqnarray}
This proves $\left(\ref{eq:612}a\right)$.  To prove $\left(\ref{eq:612}b\right)$, we note that ${\eta _j}$ is a function of the random variables $\left\{ {{Y_{j\left( {r + s} \right) + 1}},...,{Y_{j\left( {r + s} \right) + r}}} \right\}$ or ${\eta _j}$ is  $\mathcal{G}_{{u_j}}^{{v_j}}$-measurable with ${u_j} = j\left( {r + s} \right) + 1,\,{v_j} = j\left( {r + s} \right) + r$ ($\mathcal{G}_{{u_j}}^{{v_j}}$ is the $\sigma$-algebra of events generated by the random variables $\left\{ {{X_h},{U_h},u_j \le h \le v_j} \right\}$. Also ${u_{j + 1}} - {v_j} = s + 1$. Hence by Lemma 1.1 of Volkonskii and Rozannov \cite{Volkonskii1959}, we find
$$\left| {\mathbb{E}\left( {\prod\limits_{j = 0}^{k - 1} {{e^{ip{\eta _j}}}} } \right) - \prod\limits_{j = 0}^{k - 1} {\mathbb{E}\left( {{e^{ip{\eta _j}}}} \right)} } \right| \le 16k\alpha_Y \left( {s + 1} \right) \sim \frac{{16n}}{{r + s}}\alpha_Y \left( {s + 1} \right) \to 0\,\,{\rm{as}}\,\,n \to \infty, $$
the last step follows by $\left(\ref{eq:69}a\right)$ and $\left(\ref{eq:69}d\right)$.\\
Next to $\left(\ref{eq:612}c\right)$, by using stationary property of the processes $X_t, U_t$ and applying Theorem \ref{th:variance}, we obtain
$$\mathop {\inf }\limits_n \left[ {\frac{1}{n}\sum\limits_{j = 0}^{k - 1} {\mathbb{E}\left( {\eta _j^2} \right)} } \right] = \mathop {\inf }\limits_n \left[ {\frac{k}{n}\mathbb{E}\left( {\eta _0^2} \right)} \right] = \mathop {\inf }\limits_n \left[ {\frac{{k.r}}{n}\left( {{\mathop{\rm var}} \left( {{V_{n,1}}} \right) + \smallO \left( 1 \right)} \right)} \right] \to \sigma _1^2\left( x \right) \text{ as } n \to \infty. $$
Finally, to prove $\left(\ref{eq:612}d\right)$, we use
$$\left| {{U_{n,i}}} \right| \le \frac{{b_n^{{1/2} + \kappa }}}{{{b_n}}}\left\{ {\mathop {\sup }\limits_u \left| {{W_{{b_n}}}\left( u \right)} \right| + \left| {{\mu _n}} \right|} \right\} \le \frac{{b_n^\kappa }}{{b_n^{{1/2}}}}\left\{ {\frac{{\bigO\left( 1 \right)}}{{b_n^\kappa }} + \left| {{\mu _n}} \right|} \right\} \left(\text{follows by \eqref{eq:62}}\right).$$
Since ${\mu_n} \to {f_X}\left( x \right), \text{ we get } \left| {{V_{n,i}}} \right| \le \dfrac{{\bigO\left( 1 \right)}}{{b_n^{{1/2}}}}$ uniformly in $i$.\\
Hence by $\left(\ref{eq:69}c\right)$, we have $\mathop {\max }\limits_{0 \le j \le k - 1} \dfrac{{\left| {{\eta _j}} \right|}}{{\sqrt n }} \le \dfrac{{\bigO\left( 1 \right)r}}{{{{\left( {n{b_n}} \right)}^{{1/2}}}}} \to 0$ as $n \to \infty $. \\
Thus the set $\left\{ {\left| {{\eta _j}} \right| > \varepsilon {\sigma _1}\left( x \right)\sqrt n } \right\}$ becomes empty for large $n$. Therefore $\left(\ref{eq:612}d\right)$ now follows by
$$\frac{1}{n}\sum\limits_{j = 0}^{k - 1} {\mathbb{E}\left[ {\eta _j^2{I_{\left\{ {\left| {{\eta _j}} \right| > \varepsilon {\sigma _1}\left( x \right)\sqrt n } \right\}}}} \right]}  \le \frac{{\bigO\left( 1 \right)}}{n}k{\left( {\frac{r}{{b_n^{{1/2}}}}} \right)^2} \mathbb{P}\left[ {\left| {{\eta _0}} \right| > \varepsilon {\sigma _1}\left( x \right)\sqrt n } \right].$$
This completes the proof of the theorem. \hfill {\ensuremath{\square}}

\subsection{Proof of Corollary \ref{coro:asymptotic-normality}}

Note that
\begin{eqnarray}
{n^{{1/2}}}{b^{{1/2} + \kappa }}\left( {\widehat {f}_n\left(x\right) - {f_X}\left(x\right)} \right) 
&=& {n^{{1/2}}}{b^{{1/2} + \kappa }}\left[ {{{\widehat f}_n}\left(x\right) - \mathbb{E}\left( {{\widehat f}_n}\left(x\right) \right)} \right] \nonumber \\
& & \,+\,  {n^{{1/2}}}{b^{{1/2} + \kappa }}\left[ {\mathbb{E}\left( \widehat {f}_n\left(x\right) \right) - {f_X}\left(x\right)} \right]. \nonumber
\end{eqnarray}
By Theorem \ref{th:bias}, we have
\[\left| {{n^{{1/2}}}{b^{{1/2} + \kappa}}\left[ {\mathbb{E}\left({{\widehat f}_n}\left(x\right)\right) - f_X\left(x\right) } \right]} \right| = n^{1/2} b_n^{1/2 + \kappa} \bigO\left( b_n^s \right) = \bigO \left( n^{1/2} b_n^{1/2 + \kappa + s} \right) \to 0 \text{ as } n \to \infty.\]
By Theorem \ref{th: converge-a.s-estimators}, this ends the proof. \hfill {\ensuremath{\square}}

\subsection{Proof of Theorem \ref{th:upper-bounds-estimators} and \ref{th:lower-bounds-estimators}}

\noindent \textit{Proof of Theorem \ref{th:upper-bounds-estimators}}. From Theorem \ref{th:bias} and \ref{th:variance}, we have 
\[
\bias\left( {{{\widehat f }_{n}\left(x\right)}} \right) \le \bigO\left(b_n^{s} \right),
\]
and
\[
\Var \left( {{{\widehat f }_{n}\left(x\right)}} \right) \le \frac{ \bigO \left(1 \right)}{nb_n^{1 + 2\kappa}}.
\]

We choose ${b_n} = \bigO \left({n^{-1/\left(1 + 2\kappa +2s \right)}}\right)$, this yields
\begin{eqnarray}
\sup\limits_{f_X \in \mathcal{F}_{x,r}\left(A,s\right)} \mathbb{E}\left|\widehat {f}_{n}\left(x\right) - f_X\left(x\right) \right|^2 
\le \bigO \left( n^{-2s/\left(1 + 2\kappa + 2s \right)} \right).  \nonumber
\end{eqnarray}
 \hfill {\ensuremath{\square}}

\noindent \textit{Proof of Theorem \ref{th:lower-bounds-estimators}}.

Let $\mathcal{T}_{iid}$ be the set of all arbitrary estimator $\widehat{f}$ of $f_X$  based on independent identical distribution (iid) data $Y_1,...,Y_n$ from model \eqref{eq:11}. In this case, $X_1,...,X_n$ are iid and $U_1,...,U_n$ are iid. From Theorem 2 of  Belomestny and Goldenshluger \cite{Belomestny2020}, we have
\begin{eqnarray}
\inf\limits_{\widehat{f} \in \mathcal{T}_\alpha} \sup\limits_{f_X \in \mathcal{F}_{x,r}\left(A,s\right)} \mathbb{E}\left|\widehat {f}\left(x\right) - f_X\left(x\right) \right|^2 
& \ge & \inf\limits_{\widehat{f} \in \mathcal{T}_{iid}} \sup\limits_{f_X \in \mathcal{F}_{x,r}\left(A,s\right)} \mathbb{E}\left|\widehat {f}\left(x\right) - f_X\left(x\right) \right|^2   \nonumber \\
& \ge & \bigO \left( n^{-2s/\left(1 + 2\kappa + 2s \right)} \right). \nonumber
\end{eqnarray}
 \hfill {\ensuremath{\square}}

\section*{Acknowledgements}

The paper is supported by the Vietnam National Foundation for Science and Technology Development (NAFOSTED) under grand number 101.02-2023.42.

\section*{Disclosure statement}

The authors declare that they have no competing financial or non-financial interests related to this research.

\end{document}